**Unified transform method for nonautonomous linear evolution equations: rigorous analysis and qualitative theory**

Andreas Chatziafratis [1,*], Spyridon Kamvissis [2], Alain Miranville [3], Tohru Ozawa [4]

**Abstract.** Established herein is the rigorous extension of the celebrated Fokas unified transform method for linear evolution PDE with time-dependent coefficients. As illustrative examples, considered below and rigorously analysed in detail are non-homogeneous initial-boundary-value problems, formulated on the half-line, for the heat equation and the modified heat (also widely known as Barenblatt's pseudo-parabolic) equation with reaction terms and temporally varying coefficients. Solution formulae are derived formally via suitable implementation of the UTM, and rigorously verified a posteriori, in the classical sense. The latter task relies on complex-analytic tools and careful interpretation of the obtained integral representations. The explicit solutions established here are then used for investigating qualitative properties. These include regularity, boundary behavior near the axes of the domain as well as long-range and eventual periodicity asymptotics. Notably, smoothness of the solution up to the boundaries depends on compatibility conditions between initial, boundary and forcing data, into which interjected now are initial values of the coefficients. Further results pertaining to solution uniqueness, analyticity and controllability of the addressed problems are thence deduced in light of the preceding qualitative theory.

## 1. Introduction

Evolution equations with time-dependent coefficients constitute a fundamental and remarkably versatile class of mathematical models, arising naturally in the description of physical, chemical, and biological processes evolving in nonstationary environments. In many realistic settings, the mechanisms governing diffusion, reaction, transport, and relaxation are intrinsically time-dependent, reflecting changes in material properties, external forcing, environmental conditions, or the underlying dynamics of the medium. Consequently, the classical assumption of constant coefficients, although mathematically convenient, may fail to capture essential features of the phenomena under investigation. Allowing the coefficients to vary with time introduces a substantially richer mathematical structure, capable of describing transient regimes, evolving characteristic time scales, and dynamical effects that are inaccessible to autonomous models. At the same time, such generalizations pose significant analytical challenges, frequently rendering classical solution techniques ineffective and complicating the investigation of well-posedness, regularity, boundary behaviour, and long-time dynamics. These considerations provide compelling motivation for the development of systematic analytical frameworks capable of treating nonautonomous evolution equations under general initial and boundary conditions. Of particular interest are reaction–diffusion equations and Barenblatt-type pseudo-parabolic models, which provide representative yet mathematically distinct settings in which the interplay between time-dependent constitutive parameters, diffusion mechanisms, and boundary effects can be investigated in a unified manner.

From an analytical perspective, the transition from autonomous to nonautonomous evolution equations raises questions that extend well beyond the mere introduction of additional functional parameters. Indeed, the loss of temporal translation invariance fundamentally alters the structure of the associated evolution operators and, in general, precludes a direct application of the classical transform techniques that have proved so successful in the constant-coefficient setting. Particularly challenging is the construction of explicit solution representations for initial–boundary value problems, where the simultaneous presence of time-dependent coefficients and boundary constraints gives rise to intricate interactions between the temporal dynamics and the spatial spectral structure. Nevertheless, the preservation of spatial translation invariance offers an important avenue for overcoming these difficulties, suggesting that suitably adapted spectral methods may retain much of their effectiveness even in genuinely nonautonomous settings. In this direction, the unified transform method [1-10], also known as the Fokas method (see also [11-33], in alphabetical order, for a recent research programme), provides a particularly promising analytical framework, owing to its systematic exploitation of divergence-form identities, global relations, and complex spectral analysis. Extending this methodology to evolution equations with time-dependent coefficients is therefore of considerable interest, not only as a means of obtaining explicit and analytically tractable representations for specific models, but also as a step towards a broader and more flexible transform theory for linear evolution problems beyond the classical constant-coefficient paradigm.

* e-*mail address for correspondence:* chatziafrati@math.uoa.gr

[1] Waseda Research Institute for Science and Engineering, Waseda University, Tokyo, Japan
Department of Mathematics, National and Kapodistrian University of Athens, Athens, Greece

[2] Mathematics Building, University of Crete
Institute of Applied and Computational Mathematics, FORTH, Greece

[3] Henan Normal University, School of Mathematics and Statistics, Xinxiang, China
Laboratoire de Mathématiques Appliquées du Havre (LMAH), Université Le Havre Normandie, France

[4] Department of Applied Physics, Waseda University, Tokyo, Japan

The foundations of the *careful* implementation of the UTM to such problems were first laid in [13] where one of the present authors constructed explicit integral representations for solutions of nonautonomous equations (having chosen certain classes of so-called Sobolev-Galpern-type models, as examples) formulated on a semi-infinite interval, thereby properly extending the UTM to (structurally and spectrally) *admissible* linear time-dependent operators and effectively indicating the way forward for related investigations. Importantly also, substantial extensions to finite-domain problems and connections to applications in biology were successfully accomplished in [26].

The novel rigorous method is illustrated in this paper via the detailed study of two characteristic initial-boundary-value problems on the half-line. In particular:

The first one concerns the forced reaction-diffusion equation with time-dependent coefficients and is stated as follows:

***Problem 1*** *Solve*

$$(1.1)\qquad \begin{cases} \partial_t u = \alpha(t)\partial_{xx} u + \beta(t)\, u + f(x,t), \ (x,t) \in Q := \mathbb{R}^+ \times \mathbb{R}^+, \\ \lim\limits_{t\to 0^+} u(x,t) = u_0(x), \ x \in \mathbb{R}^+, \\ \lim\limits_{x\to 0^+} u(x,t) = g_0(t), \ t \in \mathbb{R}^+, \end{cases}$$

*for* $u(x,t)$.

The second problem deals with a modified heat equation, known in its simplest form [34-48] as Barenblatt equation (or two-temperature or second-order-fluid model), with time dependent coefficients, and is stated as follows:

***Problem 2*** *Solve*

$$(1.2)\qquad \begin{cases} \partial_t u = \alpha(t)\partial_{xxt} u + \beta(t)\,\partial_{xx} u + \gamma(t) u + f(x,t), \ (x,t) \in Q = \mathbb{R}^+ \times \mathbb{R}^+, \\ u(x,0) = u_0(x), \ x \in \mathbb{R}^+, \\ u(0,t) = g_0(t), \ t \in \mathbb{R}^+, \end{cases}$$

*for* $u = u(x,t)$.

***Notation and main assumptions*** Throughout this paper, we will be using the following assumptions for the given data:

(1.3) $\quad u_0(x) \in \mathcal{S}([0,\infty))$, $g_0(t) \in C^\infty([0,\infty))$, and

(1.4) $f(x,t)\in C^{\infty}(\overline{Q})$, where $Q=\mathbb{R}^{+}\times\mathbb{R}^{+}$, such that $\frac{\partial^m f}{\partial t^m}(\cdot,t)\in\mathcal{S}([0,\infty))$ ($\forall m$).

The condition imposed on $f$ means that, for every $T>0$,

$$\sup\left\{x^{\ell}\left|\frac{\partial^{n+m}f(x,t)}{\partial x^n\partial t^m}\right|:\ x>0,\ 0\le t\le T\right\}<\infty,\text{ for all nonnegative integers } n,\ m \text{ and } \ell.$$

For such data, we define the following functions:

$$\hat{u}_0(\lambda)=\int_0^{\infty}e^{-i\lambda y}u_0(y)dy \text{ and } \hat{f}(\lambda,t)=\int_0^{\infty}e^{-i\lambda y}f(y,t)dy,\text{ for } \lambda\in\mathbb{C} \text{ with } \operatorname{Im}\lambda\le 0.$$

**Some computations (1)** With $M\in\mathbb{N}$, integrating by parts, we obtain

(1.5) $$\hat{u}_0(\lambda)=\sigma_M(\lambda)+\upsilon_M(\lambda),\text{ for } \lambda\in\mathbb{C}-\{0\} \text{ with } \operatorname{Im}\lambda\le 0,$$

where we have set

$$\sigma_M(\lambda)=\sum_{k=1}^{M}\frac{u_0^{(k-1)}(0)}{(i\lambda)^k} \text{ and } \upsilon_M(\lambda)=\frac{1}{(i\lambda)^M}\int_0^{\infty}e^{-i\lambda y}u_0^{(M)}(y)dy=\frac{1}{(i\lambda)^M}\widehat{(u_0^{(M)})}(\lambda).$$

Also, we have

(1.6) $\hat{u}_0(\lambda)=\mathrm{O}(1/\lambda)$ and $\upsilon_M(\lambda)=\hat{u}_0(\lambda)-\sigma_M(\lambda)=\mathrm{O}(1/\lambda^{M+1})$, as $\lambda\to\infty$, with $\operatorname{Im}\lambda\le 0$.

**(2)** Generalizing (1.5), we see that

(1.7) $$\hat{f}(\lambda,t)=h_M(\lambda,t)+\delta_M(\lambda,t),\text{ for } \lambda\in\mathbb{C}-\{0\} \text{ with } \operatorname{Im}\lambda\le 0,$$

where

$$h_M(\lambda,t)=\sum_{k=1}^{M}\frac{1}{(i\lambda)^k}\left.\frac{\partial^{k-1}f(y,t)}{\partial y^{k-1}}\right|_{y=0} \text{ and } \delta_M(\lambda,t)\frac{1}{(i\lambda)^M}\int_0^{\infty}e^{-i\lambda y}\frac{\partial^M f(y,t)}{\partial y^M}dy.$$

Also, for every $T>0$, and uniformly for $t\le T$,

(1.8) $\hat{f}(\lambda,t)=\mathrm{O}(1/\lambda)$ and $\delta_M(\lambda,t)=\hat{f}(\lambda,t)-h_M(\lambda,t)=\mathrm{O}(1/\lambda^{M+1})$, for $\lambda\to\infty$, with $\operatorname{Im}\lambda\le 0$.

*Assumptions on the coefficients* $\alpha(t)$, $\beta(t)$, ... . For both problems, we assume that the functions $\alpha(t)$, $\beta(t)$, $\gamma(t)$ are real-valued and that

(1.9) $$\alpha(t),\ \beta(t),\ \gamma(t)\in C^{\infty}([0,\infty)).$$

In addition, we assume (for both problems) that

(1.10) $$\alpha(t)>0,\text{ for every } t\ge 0.$$

(If further assumptions are needed in order to obtain a particular result, we will state it clearly.)

Miscellaneous applications of the solutions to qualitative theory are discussed in the final section.

## 2. The heat equation: Derivation of the solution and statement of the results

Working with $\lambda\in\mathbb{C}$, $\operatorname{Im}\lambda\le 0$, and taking Fourier transforms, the differential equation in (1.1) gives

$$(2.1)\qquad \frac{\partial}{\partial t}[\hat{u}(\lambda,t)] = \alpha(t)\int_0^\infty \frac{\partial^2 u(y,t)}{\partial y^2} e^{-i\lambda y}dy + \beta(t)\int_0^\infty u(y,t)e^{-i\lambda y}dy + \int_0^\infty f(y,t)e^{-i\lambda y}dy\,,$$

since

$$\hat{u}_t = \hat{u}_t(\lambda,t) = \frac{\partial \hat{u}(\lambda,t)}{\partial t} = \int_{y=0}^{\infty} \frac{\partial u(y,t)}{\partial t} e^{-i\lambda y}dy\,.$$

Also,

$$\int_0^\infty u_x(x,t)e^{-i\lambda x}dx = -g_0(t) + i\lambda\hat{u}(\lambda,t) \text{ and } \int_0^\infty u_{xx}(x,t)e^{-i\lambda x}dx = -g_1(t) - i\lambda g_0(t) + (i\lambda)^2\hat{u}(\lambda,t)\,,$$

where we have set $g_1(t) = u_x(0,t)$. Thus, integrating by parts in (2.2), we have:

$$(2.2)\qquad \frac{\partial}{\partial t}[\hat{u}(\lambda,t)] + [\alpha(t)\lambda^2 - \beta(t)]\hat{u}(\lambda,t) = -\alpha(t)[g_1(t) + i\lambda g_0(t)] + \hat{f}(\lambda,t)\,.$$

Setting

$$\omega(\lambda,t) = \alpha(t)\lambda^2 - \beta(t) \text{ and } \Omega(\lambda,t) = \int_{\tau=0}^{t}\omega(\lambda,\tau)d\tau = \int_{\tau=0}^{t}[\alpha(\tau)\lambda^2 - \beta(\tau)]d\tau$$

and multiplying (2.2) by $e^{\Omega(\lambda,t)}$, we obtain

$$e^{\Omega(\lambda,t)}\frac{\partial}{\partial t}[\hat{u}(\lambda,t)] + e^{\Omega(\lambda,t)}[\alpha(t)\lambda^2 - \beta(t)]\hat{u}(\lambda,t) = -e^{\Omega(\lambda,t)}\alpha(t)[g_1(t) + i\lambda g_0(t)] + e^{\Omega(\lambda,t)}\hat{f}(\lambda,t)\,,$$

equivalently,

$$(2.3)\qquad \frac{\partial}{\partial t}[e^{\Omega(\lambda,t)}\hat{u}(\lambda,t)] = -e^{\Omega(\lambda,t)}\alpha(t)[g_1(t) + i\lambda g_0(t)] + e^{\Omega(\lambda,t)}\hat{f}(\lambda,t)\,.$$

Integrating (2.3), we obtain the global relation GR:

$$e^{\Omega(\lambda,t)}\hat{u}(\lambda,t) = \hat{u}_0(\lambda) - \int_{\tau=0}^{t} e^{\Omega(\lambda,\tau)}\alpha(\tau)[g_1(\tau) + i\lambda g_0(\tau)]d\tau + \int_{\tau=0}^{t} e^{\Omega(\lambda,t)}\hat{f}(\lambda,\tau)d\tau\,,$$

equivalently,

$$(2.4)\quad \hat{u}(\lambda,t) = e^{-\Omega(\lambda,t)}\hat{u}_0(\lambda) - e^{-\Omega(\lambda,t)}\int_{\tau=0}^{t} e^{\Omega(\lambda,\tau)}\alpha(\tau)[g_1(\tau) + i\lambda g_0(\tau)]d\tau + e^{-\Omega(\lambda,t)}\int_{\tau=0}^{t} e^{\Omega(\lambda,t)}\hat{f}(\lambda,\tau)d\tau\,.$$

Multiplying (2.4) by $e^{i\lambda x}$ and integrating, we find

$$(2.5)\quad 2\pi u(x,t) = \int_{-\infty}^{\infty} e^{i\lambda x}\hat{u}(\lambda,t)d\lambda = \int_{-\infty}^{\infty} e^{i\lambda x - \Omega(\lambda,t)}\hat{u}_0(\lambda)d\lambda$$

$$-\int_{-\infty}^{\infty}\left\{e^{i\lambda x - \Omega(\lambda,t)}\int_{\tau=0}^{t} e^{\Omega(\lambda,\tau)}\alpha(\tau)[g_1(\tau) + i\lambda g_0(\tau)]d\tau\right\}d\lambda + \int_{-\infty}^{\infty}\left\{e^{i\lambda x - \Omega(\lambda,t)}\int_{\tau=0}^{t} e^{\Omega(\lambda,t)}\hat{f}(\lambda,\tau)d\tau\right\}.$$

Deforming the contour of the integral which involves the terms $g_1(\tau)$ and $g_0(\tau)$, we write (2.5) as follows: For $x > 0$ and $t > 0$,

$$(2.6)\quad 2\pi u(x,t) = \int_{-\infty}^{\infty} e^{i\lambda x - \Omega(\lambda,t)}\hat{u}_0(\lambda)d\lambda$$

$$-\int_{\Gamma}\left\{e^{i\lambda x - \Omega(\lambda,t)}\int_{\tau=0}^{t} e^{\Omega(\lambda,\tau)}\alpha(\tau)[g_1(\tau) + i\lambda g_0(\tau)]d\tau\right\}d\lambda + \int_{-\infty}^{\infty}\left\{e^{i\lambda x - \Omega(\lambda,t)}\int_{\tau=0}^{t} e^{\Omega(\lambda,t)}\hat{f}(\lambda,\tau)d\tau\right\},$$

where the contour $\Gamma$ is the boundary of the set

$$\{\lambda \in \mathbb{C} : \operatorname{Im}\lambda \geq 0, \operatorname{Re}(\lambda^2) \leq 0\} = \{\lambda \in \mathbb{C} : \operatorname{Im}\lambda \geq 0, (\pi/4) \leq \arg\lambda \leq (3\pi/4),$$

as depicted in fig. 1.

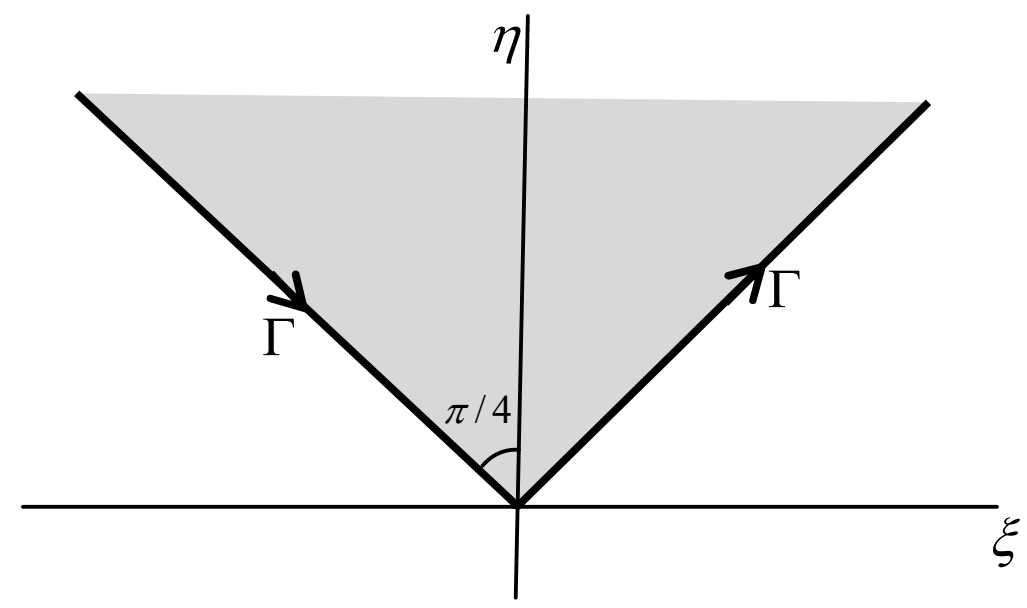

**Fig.1** *The contour* $\Gamma$ *in the* $\lambda -$*plane (* $\lambda = \xi + i\eta$ *)*

Setting $-\lambda$ in GR we find: For $\lambda \in \mathbb{C}$ with $\operatorname{Im}\lambda \geq 0$,

$$(2.7)\ \ \hat{u}(-\lambda,t) = e^{-\Omega(\lambda,t)}\hat{u}_0(-\lambda) - e^{-\Omega(\lambda,t)} \int_{\tau=0}^{t} e^{\Omega(\lambda,\tau)}\alpha(\tau)[g_1(\tau) - i\lambda g_0(\tau)]d\tau$$
$$+ e^{-\Omega(\lambda,t)} \int_{\tau=0}^{t} e^{\Omega(\lambda,t)} \hat{f}(-\lambda,\tau)d\tau .$$

This leads to the equation

$$(2.8)\ \ 0 = \int_{\Gamma} e^{i\lambda x}\hat{u}(-\lambda,t)d\lambda = \int_{\Gamma} e^{i\lambda x - \Omega(\lambda,t)}\hat{u}_0(-\lambda)d\lambda$$
$$- \int_{\Gamma} \left\{ e^{i\lambda x - \Omega(\lambda,t)} \int_{\tau=0}^{t} e^{\Omega(\lambda,\tau)}\alpha(\tau)[g_1(\tau) - i\lambda g_0(\tau)]d\tau \right\} d\lambda + \int_{\Gamma} \left\{ e^{-\Omega(\lambda,t)} \int_{\tau=0}^{t} e^{\Omega(\lambda,t)} \hat{f}(-\lambda,\tau)d\tau \right\} d\lambda .$$

Subtracting (2.8) from (2.6), we obtain: For $x > 0$ and $t > 0$,

$$(2.9)\ \ 2\pi u(x,t) = \int_{-\infty}^{\infty} e^{i\lambda x - \Omega(\lambda,t)}\hat{u}_0(\lambda)d\lambda - \int_{\Gamma} e^{i\lambda x - \Omega(\lambda,t)}\hat{u}_0(-\lambda)d\lambda$$
$$- 2\int_{\Gamma} i\lambda \left\{ e^{i\lambda x - \Omega(\lambda,t)} \int_{\tau=0}^{t} e^{\Omega(\lambda,\tau)}\alpha(\tau)g_0(\tau)d\tau \right\} d\lambda$$
$$+ \int_{-\infty}^{\infty} \left\{ e^{i\lambda x - \Omega(\lambda,t)} \int_{\tau=0}^{t} e^{\Omega(\lambda,t)} \hat{f}(\lambda,\tau)d\tau \right\} d\lambda - \int_{\Gamma} \left\{ e^{i\lambda x - \Omega(\lambda,t)} \int_{\tau=0}^{t} e^{\Omega(\lambda,t)} \hat{f}(-\lambda,\tau)d\tau \right\} d\lambda .$$

This is the solution we propose for Problem 1.

Now, fixing a $T > 0$, we can write (2.9) also in the form: For $x > 0$ and $0 < t < T$,

$$(2.10)\ \ 2\pi u(x,t) = \int_{-\infty}^{\infty} e^{i\lambda x - \Omega(\lambda,t)}\hat{u}_0(\lambda)d\lambda - \int_{\Gamma} e^{i\lambda x - \Omega(\lambda,t)}\hat{u}_0(-\lambda)d\lambda$$
$$- 2\int_{\Gamma} i\lambda \left\{ e^{i\lambda x - \Omega(\lambda,t)} \int_{\tau=0}^{T} e^{\Omega(\lambda,\tau)}\alpha(\tau)g_0(\tau)d\tau \right\} d\lambda$$
$$+ \int_{-\infty}^{\infty} \left\{ e^{i\lambda x - \Omega(\lambda,t)} \int_{\tau=0}^{t} e^{\Omega(\lambda,\tau)} \hat{f}(\lambda,\tau)d\tau \right\} d\lambda - \int_{\Gamma} \left\{ e^{i\lambda x - \Omega(\lambda,t)} \int_{\tau=0}^{t} e^{\Omega(\lambda,\tau)} \hat{f}(-\lambda,\tau)d\tau \right\} d\lambda .$$

***Notation*** We may write (2.9) as

(2.11) $$2\pi u(x,t)=[\mathcal{U}^{+}(u_0)-\mathcal{U}^{-}(u_0)-2\mathcal{G}^{t}(g_0)+\mathcal{F}^{+}(f)-\mathcal{F}^{-}(f)](x,t),$$

where

(2.12) $$\mathcal{U}^{+}(u_0),\ \mathcal{U}^{-}(u_0),\ \mathcal{G}^{t}(g_0),\ \mathcal{F}^{+}(f),\ \mathcal{F}^{-}(f),$$

are the functions of $(x,t)$, defined by the integrals in (2.9) which involve the terms

$$\hat{u}_0(\lambda),\ \hat{u}_0(-\lambda),\ g_0(\tau),\ \hat{f}(\lambda,t),\ \hat{f}(-\lambda,t),$$

respectively. Also, we may write (2.10) as

(2.13) $$2\pi u(x,t)=[\mathcal{U}^{+}(u_0)-\mathcal{U}^{-}(u_0)-2\mathcal{G}^{T}(g_0)+\mathcal{F}^{+}(f)-\mathcal{F}^{-}(f)](x,t),$$

where $\mathcal{G}^{T}(g_0)$ is the integral in (2.10) which involves the term $g_0(\tau)$.

***Proof of (2.10)*** It suffices to show that for $x>0$ and $0<t<T$,

(2.14) $$\mathcal{G}^{T}(g_0)(x,t)=\mathcal{G}^{t}(g_0)(x,t).$$

Firstly,

(2.15) $$\mathcal{G}^{T}(g_0)(x,t)-\mathcal{G}^{t}(g_0)(x,t)=-2\int_{\Gamma} i\lambda\left\{e^{i\lambda x-\Omega(\lambda,t)}\int_{\tau=t}^{T}e^{\Omega(\lambda,\tau)}\alpha(\tau)g_0(\tau)d\tau\right\}d\lambda$$

$$=-2e^{B(t)}\int_{\Gamma} i\lambda\left\{e^{i\lambda x}e^{-A(t)\lambda^2}\int_{\tau=t}^{T}e^{A(\tau)\lambda^2}e^{-B(\tau)}\alpha(\tau)g_0(\tau)d\tau\right\}d\lambda.$$

and

(2.16) $$e^{-A(t)\lambda^2}\int_{\tau=t}^{T}e^{A(\tau)\lambda^2}e^{-B(\tau)}\alpha(\tau)g_0(\tau)d\tau$$

$$=\frac{1}{\lambda^2}e^{-A(t)\lambda^2}\int_{\tau=t}^{T}\frac{d}{d\tau}[e^{A(\tau)\lambda^2}][e^{-B(\tau)}g_0(\tau)]d\tau=\mathrm{O}(1/\lambda^2),\text{ as }\lambda\to\infty\text{ with }\frac{\pi}{4}\le\arg\lambda\le\frac{3\pi}{4},$$

as integration by parts shows.

It follows from (2.16) that

$$\lim_{R\to\infty}\int_{\substack{|\lambda|=R\\(\pi/4)\le\arg\lambda\le(3\pi/4)}} i\lambda\left\{e^{i\lambda x}e^{-A(t)\lambda^2}\int_{\tau=t}^{T}e^{A(\tau)\lambda^2}e^{-B(\tau)}\alpha(\tau)g_0(\tau)d\tau\right\}d\lambda=0.$$

Therefore, (2.14) follows from (2.15) and Cauchy's theorem.

We will verify that, indeed, the function $u(x,t)$, which we derived above, solves Problem 1, by proving the following theorem. We recall the main assumptions: (1.3), (1.4), (1.9) and (1.10).

**Theorem 1** *The function* $u(x,t)$*, defined, for* $(x,t)\in\mathbb{R}^{+}\times\mathbb{R}^{+}$*, by (2.9) – equivalently by (2.10) – solves Problem 1. Moreover we have:*

*1*$^{st}$ $u(x,t)$ *is* $C^{\infty}$ *for* $(x,t)\in\mathbb{R}^{+}\times\mathbb{R}^{+}$.

*2*$^{nd}$ *The convergence* $\lim_{t\to 0^{+}}u(x,t)=u_0(x)$ *is uniform for* $x$ *in compact subsets of* $\mathbb{R}^{+}$.

*3rd The convergence* $\lim_{x\to 0^+} u(x,t) = g_0(t)$ *is uniform for* $t$ *in compact subsets of* $\mathbb{R}^+$.

Further properties of the solution are given by the following theorems. These properties describe the boundary behaviour of the solution as well as its asymptotic behaviour as $x \to +\infty$.

**Theorem 2** *For the function* $u(x,t)$, *defined, for* $(x,t) \in \mathbb{R}^+ \times \mathbb{R}^+$, *by (2.9), we have:*

*1st The limits*

$$U_{n,m}(x) := \lim_{t\to 0^+} \frac{\partial^{n+m} u(x,t)}{\partial x^n \partial t^m}, \quad n, m = 0,1,2, ...,$$

*exist and are uniform for* $x$ *in compact subsets of* $\mathbb{R}^+$. *Moreover the functions* $U_{n,m}(x)$ *are* $C^\infty$ *for* $x \in [0,\infty)$.

*2nd The limits*

$$G_{n,m}(t) := \lim_{x\to 0^+} \frac{\partial^{n+m} u(x,t)}{\partial x^n \partial t^m}, \quad n, m = 0,1,2, ...,$$

*exist and are uniform for* $t$ *in compact subsets of* $\mathbb{R}^+$. *Moreover the functions* $G_{n,m}(t)$ *are* $C^\infty$ *for* $t \in [0,\infty)$.

*3rd* $u(x,t)$, *extended to* $\overline{Q} - \{(0,0)\}$ *by setting*

$$u^*(x,0) := u_0(x) \text{ for } x > 0 \text{ and } u^*(0,t) := g_0(t) \text{ for } t > 0,$$

*is* $C^\infty$ *in* $\overline{Q} - \{(0,0)\}$.

*4th* $\lim_{t\to 0^+} \frac{\partial^n u(x,t)}{\partial x^n} = \frac{d^n u_0(x)}{dx^n}$, *for* $x > 0$, *and* $\lim_{x\to 0^+} \frac{\partial^n u(x,t)}{\partial t^n} = \frac{d^n g_0(x)}{dt^n}$, *for* $t > 0$.

The following theorem gives conditions for the existence of the limits of the functions

$$u(x,t), u_x(x,t), u_{xx}(x,t), u_{xxx}(x,t), u_t(x,t),$$

as $(x,t) \to (0^+, 0^+)$. The existence of the corresponding limits for higher order derivatives of $u(x,t)$ require further compatibility conditions for the data at the point $(0,0)$.

**Theorem 3** *The functions* $u(x,t)$, $\partial u(x,t)/\partial x$, $\partial^2 u(x,t)/\partial x^2$ *and* $\partial^3 u(x,t)/\partial x^3$ *(originally defined by (2.9) for* $(x,t) \in Q := \mathbb{R}^+ \times \mathbb{R}^+$ *and extended to* $\overline{Q} - \{(0,0)\}$ *by Theorem 2) satisfy the following:*

*1st If* $u_0(0) = g_0(0)$ *then*

$$\lim_{\overline{Q} \ni (x,t) \to (0,0)} u(x,t) = u_0(0) \text{ and } \lim_{\overline{Q} \ni (x,t) \to (0,0)} \frac{\partial u(x,t)}{\partial x} = \left.\frac{du_0(x)}{dx}\right|_{x=0}.$$

*2nd If*

$$u_0(0) = g_0(0) \text{ and } \left.\frac{dg_0(t)}{dt}\right|_{t=0} = \alpha(0) \left.\frac{d^2 u_0(x)}{dx^2}\right|_{x=0} + \beta(0) g_0(0) + f(0,0) \qquad (2.17)$$

*then*

$$\lim_{\overline{Q}\ni(x,t)\to(0,0)} \frac{\partial^2 u(x,t)}{\partial x^2} = \left.\frac{d^2 u_0(x)}{dx^2}\right|_{x=0} \quad \textit{and} \quad \lim_{\overline{Q}\ni(x,t)\to(0,0)} \frac{\partial^3 u(x,t)}{\partial x^3} = \left.\frac{d^3 u_0(x)}{dx^3}\right|_{x=0}.$$

*3rd Under the assumption (2.17) of 2nd part,*

$$\lim_{\overline{Q}\ni(x,t)\to(0,0)} \frac{\partial u(x,t)}{\partial t} = \alpha(0)\left.\frac{d^2 u_0(x)}{dx^2}\right|_{x=0} + \beta(0)u_0(0) + f(0,0).$$

**Theorem 4** *For the function* $u(x,t)$*, defined, for* $(x,t)\in\mathbb{R}^+\times\mathbb{R}^+$*, by (2.9), we have: For* $T>0$,

$$\lim_{x\to+\infty}\left( x^{\ell}\,\frac{\partial^{n+m}u(x,t)}{\partial x^n \partial t^m}\right) = 0,$$

*for nonnegative integers* $n$, $m$ *and* $\ell$*, uniformly for* $0<t\le T$.

## 3. Proofs of Theorems 1 - 4

**Proof of Theorem 1** The theorem follows from the conclusions of the following steps combined with the specific computations which are carried out in their proofs.

***Step 1*** *The functions defined by the integrals (2.12) are* $C^\infty$ *for* $(x,t)\in\mathbb{R}^+\times\mathbb{R}^+$.

*Proof* That the integrals $\mathcal{U}^-(u_0)$, $\mathcal{G}^t(g_0)$, $\mathcal{F}^-(f)$ define $C^\infty$ functions for $(x,t)\in\mathbb{R}^+\times\mathbb{R}^+$, it follows from the presence in these integrals of the factor $e^{i\lambda x}$ (as long as $x$ is kept away from $0$), and the fact that these integrals are taken for $\lambda\in\Gamma$, since, for these $\lambda's$,

$$\left|e^{i\lambda x}\right| = e^{-(\operatorname{Im}\lambda)x} = e^{-x(|\lambda|/\sqrt{2})}.$$

To deal with the integral $\mathcal{U}^+(u_0)$, it suffices, using (1.6), to write it as follows

$$\text{(3.1)}\quad \mathcal{U}^+(u_0)(x,t) = \int_{-1}^{1} e^{i\lambda x-\Omega(\lambda,t)}\hat{u}_0(\lambda)d\lambda$$
$$+\sum_{k=1}^{M} u_0^{(k-1)}(0)\left(\int_{-\infty}^{-1}+\int_{1}^{\infty}\right)e^{i\lambda x-\Omega(\lambda,t)}\frac{d\lambda}{(i\lambda)^k} + \left(\int_{-\infty}^{-1}+\int_{1}^{\infty}\right)e^{i\lambda x-\Omega(\lambda,t)}(u_0^{(M)})\hat{}(\lambda)\frac{d\lambda}{(i\lambda)^M} \quad (\forall M)$$

and to observe that

$$\text{(3.2)}\qquad \left(\int_{-\infty}^{-1}+\int_{1}^{\infty}\right)e^{i\lambda x-\omega(\lambda)t}\frac{d\lambda}{(i\lambda)^k} = \int_{\Gamma_1} e^{i\lambda x-\omega(\lambda)t}\frac{d\lambda}{(i\lambda)^k} + \int_{\varepsilon_1} e^{i\lambda x-\omega(\lambda)t}\frac{d\lambda}{(i\lambda)^k} + \int_{\varepsilon_2} e^{i\lambda x-\omega(\lambda)t}\frac{d\lambda}{(i\lambda)^k}.$$

where $\Gamma_1$ is the part of $\Gamma$ where $|\lambda|\ge 1$, and $\varepsilon_1$, $\varepsilon_2$ are the straight line segments

$$\varepsilon_1 = [(-1/\sqrt{2})+(i/\sqrt{2}), -1/\sqrt{2}] \text{ and } \varepsilon_2 = [1/\sqrt{2}, (1/\sqrt{2})+(i/\sqrt{2})]. \text{ (See fig. 2.)}$$

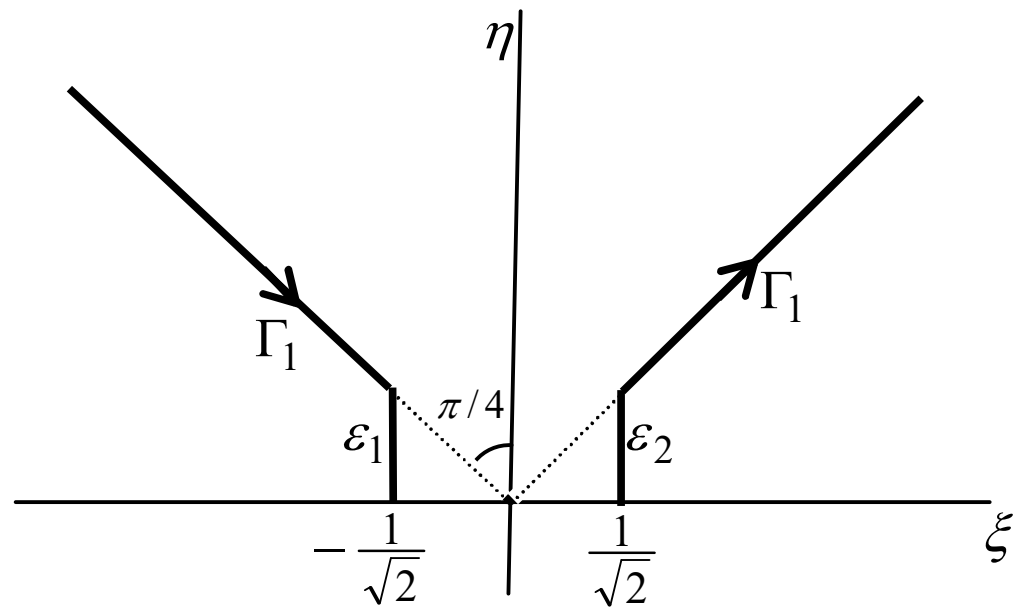

**Fig.2** *The two pieces of the contour* $\Gamma_1$ *and the line segments* $\varepsilon_1$, $\varepsilon_2$

Indeed, (3.1) and (3.2) imply that, for sufficiently large $M$,

$$(3.3)\quad \frac{\partial^{n+m}\mathcal{U}^+(u_0)(x,t)}{\partial x^n \partial t^m} = \int_{-1}^{1} (i\lambda)^n e^{i\lambda x} \frac{\partial^m [e^{-\Omega(\lambda,t)}]}{\partial t^m} \hat{u}_0(\lambda) d\lambda$$

$$+ \sum_{k=1}^{M} u_0^{(k-1)}(0) \left( \int_{\Gamma_1} + \int_{\varepsilon_1+\varepsilon_2} \right) (i\lambda)^n e^{i\lambda x} \frac{\partial^m [e^{-\Omega(\lambda,t)}]}{\partial t^m} \frac{d\lambda}{(i\lambda)^k}$$

$$+ \left( \int_{-\infty}^{-1} + \int_{1}^{\infty} \right) (i\lambda)^n e^{i\lambda x} \frac{\partial^m [e^{-\Omega(\lambda,t)}]}{\partial t^m} (u_0^{(M)})\hat{}(\lambda) \frac{d\lambda}{(i\lambda)^M} .$$

We emphasise that, in (3.1), (3.2) and (3.3), we assume $x > 0$ and $t > 0$.

To deal with the integral $\mathcal{F}^+(f)$, we use the formula $\hat{f}(\lambda,t) = h_M(\lambda,t) + \delta_M(\lambda,t)$ of (1.7), and an analogue of (3.1), and we write

$$(3.4)\quad (\mathcal{F}^+(f))(x,t) = \int_{-1}^{1} e^{i\lambda x - \Omega(\lambda,t)} \tilde{\hat{f}}(\lambda,\Omega(\lambda,\cdot),t) d\lambda$$

$$+ \left( \int_{\Gamma_1} + \int_{\varepsilon_1+\varepsilon_2} \right) e^{i\lambda x - \Omega(\lambda,t)} \tilde{h}_M(\lambda,\Omega(\lambda,\cdot),t) d\lambda + \left( \int_{-\infty}^{-1} + \int_{1}^{\infty} \right) e^{i\lambda x - \Omega(\lambda,t)} \tilde{\delta}_M(\lambda,\Omega(\lambda,\cdot),t) d\lambda ,$$

where

$$\tilde{\hat{f}}(\lambda,\Omega(\lambda,\cdot),t) = \int_{\tau=0}^{t} e^{\Omega(\lambda,\tau)} \hat{f}(\lambda,\tau) d\tau ,$$

$$\tilde{h}_M(\lambda,\Omega(\lambda,\cdot),t) = \int_{\tau=0}^{t} e^{\Omega(\lambda,\tau)} h_M(\lambda,\tau) d\tau ,$$

$$\tilde{\delta}_M(\lambda,\Omega(\lambda,\cdot),t) = \int_{\tau=0}^{t} e^{\Omega(\lambda,\tau)} \delta_M(\lambda,\tau) d\tau .$$

Then, for sufficiently large $M$, (3.4) leads to

$$(3.5)\quad \frac{\partial^{n+m}\mathcal{F}^+(f)(x,t)}{\partial x^n \partial t^m} = \int_{-1}^{1} (i\lambda)^n e^{i\lambda x} \frac{\partial^m [e^{-\Omega(\lambda,t)} \tilde{\hat{f}}(\lambda,\Omega(\lambda,\cdot),t)]}{\partial t^m} d\lambda$$

$$+ \left( \int_{\Gamma_1} + \int_{\varepsilon_1+\varepsilon_2} \right) (i\lambda)^n e^{i\lambda x} \frac{\partial^m [e^{-\Omega(\lambda,t)} \tilde{h}_M(\lambda,\Omega(\lambda,\cdot),t)]}{\partial t^m} d\lambda$$

$$+\left(\int_{-\infty}^{-1}+\int_{1}^{\infty}\right)(i\lambda)^n e^{i\lambda x}\frac{\partial^m[e^{-\Omega(\lambda,t)}\tilde{\delta}_M(\lambda,\Omega(\lambda,\cdot),t)]}{\partial t^m}d\lambda .$$

As we have proved, the quantities $\mathcal{G}^t(g_0)(x,t)$ and $\mathcal{G}^T(g_0)(x,t)$ are equal for $x>0$ and $0<t<T$. However, their evaluations for $x=0$ are, in general, different. More precisely, we have

***Step 2*** *The integrals* $\mathcal{G}^t(g_0)(x,t)\Big|_{x=0}$ *and* $\mathcal{G}^T(g_0)(x,t)\Big|_{x=0}$ *exist in the generalized sense and*

(3.6) $$\mathcal{G}^T(g_0)(x,t)\Big|_{x=0}=-\pi g_0(t),\ \textit{for}\ 0<t<T,$$

*while*

(3.7) $$\mathcal{G}^t(g_0)(x,t)\Big|_{x=0}=-\frac{1}{2}\pi g_0(t),\ \textit{for}\ t>0.$$

*Proof* Let us set

$$\mu=i\lambda^2 \text{ and } \sigma=\varphi(\tau)=\int_{\varsigma=0}^{\tau}\alpha(\varsigma)d\varsigma .$$

Then $\tau=\psi(\sigma)$, where $\psi:=\varphi^{-1}$,

$$d\mu=2i\lambda d\lambda \text{ and } d\sigma=\alpha(\tau)d\tau .$$

Thus, for a continuous function $h$ of $\tau$,

(3.8) $$I_{h,T}(t):=\int_\Gamma i\lambda\left\{e^{-\lambda^2\varphi(t)}\int_{\tau=0}^{T}e^{\lambda^2\varphi(\tau)}h(\tau)d\tau\right\}d\lambda=-\frac{1}{2}\int_{\mu=-\infty}^{\infty}\left\{e^{i\mu\varphi(t)}\int_{\tau=0}^{T}e^{-i\mu\varphi(\tau)}h(\tau)d\tau\right\}d\mu$$

and

$$\int_{\tau=0}^{T}e^{-i\mu\varphi(\tau)}h(\tau)d\tau=\int_{\sigma=0}^{\varphi(T)}e^{-i\mu\sigma}h(\psi(\sigma))\frac{d\sigma}{\alpha(\psi(\sigma))}=\hat{H}_T(\mu),$$

where

$$H_T(\sigma)=\begin{cases}\dfrac{h(\psi(\sigma))}{\alpha(\psi(\sigma))} & \textit{for}\ \sigma\in[0,\varphi(T)]\\ 0 & \textit{for}\ \sigma\notin[0,\varphi(T)].\end{cases}$$

Since $0<\varphi(t)<\varphi(T)$, it follows, from (3.8) and Fourier's inversion formula, that

(3.9) $$I_{h,T}(t)=-\frac{1}{2}2\pi H_T(\varphi(t))=-\pi\frac{h(t)}{\alpha(t)},\ \text{for}\ 0<t<T.$$

Now we prove (3.6). We have

(3.10) $$\mathcal{G}^T(g_0)(x,t)\Big|_{x=0}=\int_\Gamma i\lambda\left\{e^{-\Omega(\lambda,t)}\int_{\tau=0}^{T}e^{\Omega(\lambda,\tau)}\alpha(\tau)g_0(\tau)d\tau\right\}d\lambda$$

$$=\exp\left(\int_{\tau=0}^{t}\beta(\tau)d\tau\right)\int_\Gamma i\lambda\left\{e^{-\lambda^2\varphi(t)}\int_{\tau=0}^{T}e^{\lambda^2\varphi(\tau)}\left[\exp\left(-\int_{\varsigma=0}^{\tau}\beta(\varsigma)d\varsigma\right)\alpha(\tau)g_0(\tau)\right]d\tau\right\}d\lambda$$

$$= \exp\left( \int_{\tau=0}^{t} \beta(\tau) d\tau \right) I_{h,T}(t),$$

when

$$h(\tau) = \exp\left( - \int_{\varsigma=0}^{\tau} \beta(\varsigma) d\varsigma \right) \alpha(\tau) g_0(\tau).$$

Thus, (3.6) follows from (3.9) and (3.10).

The proof of (3.7) is similar. More precisely, we can show that

$$\text{(3.11)} \qquad I_{h,t}(t) := \int_{\Gamma} i\lambda \left\{ e^{-\lambda^2 \varphi(t)} \int_{\tau=0}^{t} e^{\lambda^2 \varphi(\tau)} h(\tau) d\tau \right\} d\lambda = -\frac{1}{2}\pi \frac{h(t)}{\alpha(t)}.$$

The factor $1/2$ in (3.11) is due to the fact that in this case we apply Fourier's inversion formula to the function

$$H_t(\sigma) = \begin{cases} \dfrac{h(\psi(\sigma))}{\alpha(\psi(\sigma))} & \text{for } \sigma \in [0, \varphi(t)] \\ 0 & \text{for } \sigma \notin [0, \varphi(t)] \end{cases}$$

and

$$\int_{\mu=-\infty}^{\infty} e^{i\mu\varphi(t)} \hat{H}_t(\mu) d\mu = 2\pi \frac{1}{2} H_t(\varphi(t)).$$

Finally, we point out that the existence of the integrals $\mathcal{G}^t(g_0)(x,t)\big|_{x=0}$ and $\mathcal{G}^T(g_0)(x,t)\big|_{x=0}$, in the generalized sense, is part of Fourier's inversion formula. (It is because of this fact that we can set $t = 0$ in the integrals $\mathcal{G}^t(g_0)(x,t)$ and $\mathcal{G}^T(g_0)(x,t)$.)

***Step 3*** *For* $x > 0$,

$$\text{(3.12)} \qquad \lim_{t \to 0^+} \mathcal{U}^+(u_0)(x,t) = 2\pi u_0(x) \text{ and } \lim_{t \to 0^+} \mathcal{U}^-(u_0)(x,t) = 0.$$

*Proof* Since the evaluation

$$\text{(3.13)} \qquad \mathcal{U}^+(u_0)(x,t)\big|_{t=0} = \int_{-\infty}^{\infty} e^{i\lambda x} \hat{u}_0(\lambda) d\lambda = 2\pi u_0(x),$$

it suffices to justify, in (3.12), the interchange of the limit with the integral defining $\mathcal{U}^+(u_0)$. For this let us recall (see (1.7)) that

$$\hat{u}_0(\lambda) = \frac{u_0(0)}{i\lambda} + \frac{(u_0')\hat{}(\lambda)}{i\lambda},$$

write

$$\text{(3.14)} \quad \mathcal{U}^+(u_0)(x,t) = \int_{-1}^{1} e^{i\lambda x - \Omega(\lambda,t)} \hat{u}_0(\lambda) d\lambda + \left( \int_{-\infty}^{-1} + \int_{1}^{\infty} \right) e^{i\lambda x - \Omega(\lambda,t)} \frac{u_0(0)}{i\lambda} d\lambda$$
$$+ \left( \int_{-\infty}^{-1} + \int_{1}^{\infty} \right) e^{i\lambda x - \Omega(\lambda,t)} \frac{(u_0')\hat{}(\lambda)}{i\lambda} d\lambda$$

and observe that

$$\text{(3.15)} \quad \lim_{t \to 0^+} \left( \int_{-\infty}^{-1} + \int_{1}^{\infty} \right) e^{i\lambda x - \Omega(\lambda,t)} \frac{u_0(0)}{i\lambda} d\lambda = \lim_{t \to 0^+} \int_{\Gamma_1 + \varepsilon_1 + \varepsilon_2} e^{i\lambda x - \Omega(\lambda,t)} \frac{u_0(0)}{i\lambda} d\lambda$$

$$= \int_{\Gamma_1+\varepsilon_1+\varepsilon_2} e^{i\lambda x}\frac{u_0(0)}{i\lambda}d\lambda = \left(\int_{-\infty}^{-1} + \int_1^{\infty}\right) e^{i\lambda x}\frac{u_0(0)}{i\lambda}d\lambda,$$

where we used also Jordan's lemma.

Thus, the first equation in (3.12) follows from (3.13), (3.14) and (3.15). The second one follows immediately from the fact that the integration is taken on $\Gamma$ and that the presence of the factor $e^{i\lambda x}$ implies that

$$\lim_{t\to 0^+} \mathcal{U}^-(u_0)(x,t) = \int_\Gamma e^{i\lambda x}\hat{u}_0(-\lambda)d\lambda = 0.$$

***Step 4*** *We have*

(3.16) $$\lim_{x\to 0^+} \mathcal{G}^T(g_0)(x,t) = \mathcal{G}^T(g_0)(x,t)\Big|_{x=0}.$$

*Proof* Integration by parts shows that, for $\lambda \neq 0$,

$$\int_{\tau=0}^{T} e^{\Omega(\lambda,\tau)}\alpha(\tau)g_0(\tau)d\tau = \frac{1}{\lambda^2}\int_{\tau=0}^{T}\frac{d}{d\tau}[e^{A(\tau)\lambda^2}]e^{-B(\tau)}g_0(\tau)d\tau$$

$$= \frac{1}{\lambda^2}e^{\Omega(\lambda,T)}g_0(T) - \frac{1}{\lambda^2}g_0(0) - \frac{1}{\lambda^2}\int_{\tau=0}^{T} e^{A(\tau)\lambda^2}\frac{d}{d\tau}[e^{-B(\tau)}g_0(\tau)]d\tau,$$

where we have set

$$A(t) = \int_{\tau=0}^{t}\alpha(\tau)d\tau \text{ and } B(t) = \int_{\tau=0}^{t}\beta(\tau)d\tau.$$

Therefore, for $\lambda \neq 0$,

$$i\lambda\left\{e^{i\lambda x-\Omega(\lambda,t)}\int_{\tau=0}^{T} e^{\Omega(\lambda,\tau)}\alpha(\tau)g_0(\tau)d\tau\right\}$$

$$= \frac{i}{\lambda}e^{i\lambda x}e^{\Omega(\lambda,T)-\Omega(\lambda,t)}g_0(T) - \frac{i}{\lambda}e^{i\lambda x-\Omega(\lambda,t)}g_0(0) - \frac{i}{\lambda}e^{i\lambda x-\Omega(\lambda,t)}\int_{\tau=0}^{T} e^{A(\tau)\lambda^2}\frac{d}{d\tau}[e^{-B(\tau)}g_0(\tau)]d\tau.$$

Thus, writing

$$\mathcal{G}^T(g_0)(x,t) = \int_\Gamma i\lambda\left\{e^{i\lambda x-\Omega(\lambda,t)}\int_{\tau=0}^{T} e^{\Omega(\lambda,\tau)}\alpha(\tau)g_0(\tau)d\tau\right\}d\lambda = \int_{\Gamma_0}\cdots d\lambda + \int_{\Gamma_1}\cdots d\lambda,$$

(where $\Gamma_0 := \Gamma \cap \{|\lambda| \leq 1\}$) we see that, in order to prove (3.16), it suffices to show the following equations:

(3.17) $$\lim_{x\to 0^+}\int_{\Gamma_1}\frac{i}{\lambda}e^{i\lambda x}e^{\Omega(\lambda,T)-\Omega(\lambda,t)}g_0(T)d\lambda = \int_{\Gamma_1}\frac{i}{\lambda}e^{\Omega(\lambda,T)-\Omega(\lambda,t)}g_0(T)d\lambda$$

(3.18) $$\lim_{x\to 0^+}\int_{\Gamma_1}\frac{i}{\lambda}e^{i\lambda x}e^{-\Omega(\lambda,t)}g_0(0)d\lambda = \int_{\Gamma_1}\frac{i}{\lambda}e^{-\Omega(\lambda,t)}g_0(0)d\lambda$$

(3.19) $$\lim_{x\to 0^+}\int_{\Gamma_1}\frac{i}{\lambda}e^{i\lambda x-\Omega(\lambda,t)}\int_{\tau=0}^{T} e^{A(\tau)\lambda^2}\frac{d}{d\tau}[e^{-B(\tau)}g_0(\tau)]d\tau = \int_{\Gamma_1}\frac{i}{\lambda}e^{-\Omega(\lambda,t)}\int_{\tau=0}^{T} e^{A(\tau)\lambda^2}\frac{d}{d\tau}[e^{-B(\tau)}g_0(\tau)]d\tau$$

We point out that, as it will become clear from the proofs of the above equations, the integrals in the RHSs of (3.17) and (3.18) exist in the generalised sense, while the integral in the RHS of (3.19) converges absolutely.

Firstly, integrating by parts once more, we see that

(3.20) $$\frac{i}{\lambda}e^{i\lambda x-\Omega(\lambda,t)}\int_{\tau=0}^{T} e^{A(\tau)\lambda^2}\frac{d}{d\tau}[e^{-B(\tau)}g_0(\tau)]d\tau = \mathrm{O}(1/\lambda^3)\text{, as } \lambda\to\infty \text{ with } \lambda\in\Gamma_1,$$

since, for these $\lambda's$, $\left|e^{A(t)\lambda^2}\right|=1,\ \forall t$. This implies immediately (3.19).

*Proof of (3.17)* It suffices to show that

(3.21) $$\lim_{x\to 0^+}\int_{\Gamma_1}\frac{i}{\lambda}e^{i\lambda x}e^{\lambda^2[A(T)-A(t)]}d\lambda = \int_{\Gamma_1}\frac{i}{\lambda}e^{\lambda^2[A(T)-A(t)]}d\lambda.$$

Since $\left|e^{\lambda^2[A(T)-A(t)]}\right|\le 1$, for $(\pi/4)\le \arg\lambda\le(3\pi/4)$,

$$\int_{\Gamma_1}\frac{i}{\lambda}e^{i\lambda x}e^{\lambda^2[A(T)-A(t)]}d\lambda = \int_{\substack{|\lambda|=1\\(\pi/4)\le\arg\lambda\le(3\pi/4)}}\frac{i}{\lambda}e^{i\lambda x}e^{\lambda^2[A(T)-A(t)]}d\lambda,$$

hence

(3.22) $$\lim_{x\to 0^+}\int_{\Gamma_1}\frac{i}{\lambda}e^{i\lambda x}e^{\lambda^2[A(T)-A(t)]}d\lambda = \int_{\substack{|\lambda|=1\\(\pi/4)\le\arg\lambda\le(3\pi/4)}}\frac{i}{\lambda}e^{\lambda^2[A(T)-A(t)]}d\lambda.$$

On the other hand, setting $\mu=-i\lambda^2$, we find

(3.23) $$\int_{\substack{|\lambda|=1\\(\pi/4)\le\arg\lambda\le(3\pi/4)}}\frac{i}{\lambda}e^{\lambda^2[A(T)-A(t)]}d\lambda = -\frac{1}{2}\int_{\substack{|\mu|=1\\0\le\arg\mu\le\pi}}\frac{1}{i\mu}e^{i\mu[A(T)-A(t)]}d\mu$$
$$= -\frac{1}{2}\int_{\substack{\mu\in\mathbb{R}\\|\mu|\ge 1}}\frac{1}{i\mu}e^{i\mu[A(T)-A(t)]}d\mu = \int_{\Gamma_1}\frac{i}{\lambda}e^{\lambda^2[A(T)-A(t)]}d\lambda.$$

(We point out that the above integrals exist in the generalized sense.)
Thus, (3.21) follows from (3.22) and (3.23).

*Proof of (3.18)* It suffices to show that

(3.24) $$\lim_{x\to 0^+}\int_{\Gamma_1}\frac{i}{\lambda}e^{i\lambda x}e^{-\lambda^2 A(t)}d\lambda = \int_{\Gamma_1}\frac{i}{\lambda}e^{-\lambda^2 A(t)}d\lambda.$$

Since $\left|e^{-\lambda^2 A(t)}\right|\le 1$, for $0\le\arg\lambda\le(\pi/4)$ and for $(3\pi/4)\le\arg\lambda\le\pi$,

$$\int_{\Gamma_1}\frac{i}{\lambda}e^{i\lambda x}e^{-\lambda^2 A(t)}d\lambda = \left(\int_{\substack{\lambda\in\mathbb{R}\\|\lambda|\ge 1}} - \int_{\substack{|\lambda|=1\\(3\pi/4)\le\arg\lambda\le\pi}} - \int_{\substack{|\lambda|=1\\0\le\arg\lambda\le(\pi/4)}}\right)\left[\frac{i}{\lambda}e^{i\lambda x}e^{-\lambda^2 A(t)}d\lambda\right],$$

and, therefore,

$$\lim_{x\to 0^+}\int_{\Gamma_1}\frac{i}{\lambda}e^{i\lambda x}e^{-\lambda^2 A(t)}d\lambda = \left(\int_{\substack{\lambda\in\mathbb{R}\\|\lambda|\ge 1}} - \int_{\substack{|\lambda|=1\\(3\pi/4)\le\arg\lambda\le\pi}} - \int_{\substack{|\lambda|=1\\0\le\arg\lambda\le(\pi/4)}}\right)\left[\frac{i}{\lambda}e^{-\lambda^2 A(t)}d\lambda\right]$$

$$= -\left( \int\limits_{\substack{|\lambda|=1 \\ (3\pi/4)\le \arg\lambda \le \pi}} + \int\limits_{\substack{|\lambda|=1 \\ 0\le \arg\lambda \le (\pi/4)}} \right) \left[ \frac{i}{\lambda} e^{-\lambda^2 A(t)} d\lambda \right]$$

$$= -\frac{1}{2}\left( \int\limits_{\substack{|\mu|=1 \\ 0\le \arg\mu \le (\pi/2)}} + \int\limits_{\substack{|\mu|=1 \\ (\pi/2)\le \arg\mu \le \pi}} \right) \left[ \frac{1}{-i\mu} e^{i\mu A(t)} d\mu \right]$$

$$= -\frac{1}{2}\left( \int\limits_{\mu=-\infty}^{-1} + \int\limits_{\mu=1}^{\infty} \right) \left[ \frac{1}{-i\mu} e^{i\mu A(t)} d\mu \right] = \int\limits_{\Gamma_1} \frac{i}{\lambda} e^{-\lambda^2 A(t)} d\lambda ,$$

where we also set $\mu = i\lambda^2$. This proves (3.24) and completes the proof of (3.18).

**Proof of Theorem 2** ***Step 1*** The 1st assertion follows from (3.3) and (3.5) as far as the integrals $\mathcal{U}^{+}(u_0)$ and $\mathcal{F}^{+}(f)$ are concerned. That this holds also for the integrals $\mathcal{U}^{-}(u_0)$, $\mathcal{G}^{t}(g_0)$, $\mathcal{F}^{-}(f)$, follows from the presence in these integrals of the factor $e^{i\lambda x}$.

***Step 2*** The 2nd assertion for the integral $\mathcal{U}^{+}(u_0)$ follows from the fact that this integral is taken for $\lambda \in \mathbb{R}$, combined with the presence in this integral of the factor

$$e^{i\lambda x - \Omega(\lambda,t)} = e^{i\lambda x} e^{B(t)} e^{-A(t)\lambda^2} ,$$

as long as $t$ is kept away from $0$. Similarly we can deal with the integral $\mathcal{U}^{-}(u_0)$, since

$$\mathcal{U}^{-}(u_0)(x,t) = \int\limits_{\Gamma} e^{i\lambda x - \Omega(\lambda,t)} \hat{u}_0(-\lambda) d\lambda = \int\limits_{-\infty}^{\infty} e^{i\lambda x - \Omega(\lambda,t)} \hat{u}_0(-\lambda) d\lambda .$$

For the integrals $\mathcal{G}^{t}(g_0)$, $\mathcal{F}^{+}(f)$ and $\mathcal{F}^{-}(f)$, however, we need to write them in another form, which we describe in the following steps.

***Step 3*** Integration by parts shows

$$e^{-A(t)\lambda^2} \int\limits_{\tau=0}^{t} e^{A(\tau)\lambda^2} e^{-B(\tau)} \alpha(\tau) g_0(\tau) d\tau = \frac{1}{\lambda^2} e^{-A(t)\lambda^2} \int\limits_{\tau=0}^{t} \frac{d}{d\tau}[e^{A(\tau)\lambda^2}][e^{-B(\tau)} g_0(\tau)] d\tau$$

$$= \frac{G_0(t)}{\lambda^2} - \frac{g_0(0)}{\lambda^2} e^{-A(t)\lambda^2} - \frac{1}{\lambda^2} e^{-A(t)\lambda^2} \int\limits_{\tau=0}^{t} e^{A(\tau)\lambda^2} \frac{d}{d\tau}[G_0(t)] d\tau ,$$

where $G_0(t) = e^{-B(t)} g_0(t)$.
Further integration by parts gives

$$(3.25)\quad e^{-A(t)\lambda^2} \int\limits_{\tau=0}^{t} e^{A(\tau)\lambda^2} e^{-B(\tau)} \alpha(\tau) g_0(\tau) d\tau = \sum_{l=1}^{L} (-1)^{l-1} \left[ \frac{1}{\lambda^{2l}} G_{l-1}(t) - \frac{1}{\lambda^{2l}} G_{l-1}(0) e^{-A(t)\lambda^2} \right]$$

$$+ (-1)^{L} \frac{1}{\lambda^{2L}} e^{-A(t)\lambda^2} \int\limits_{\tau=0}^{t} e^{A(\tau)\lambda^2} \frac{d}{d\tau}[G_{L-1}(t)] d\tau ,$$

where

$$G_1(t) = \frac{1}{\alpha(t)} \frac{d}{dt}[e^{-B(t)} g_0(t)] = \frac{1}{\alpha(t)} \frac{d}{dt}[G_0(t)],\ G_2(t) = \frac{1}{\alpha(t)} \frac{d}{dt}[G_1(t)],\ G_3(t) = \frac{1}{\alpha(t)} \frac{d}{dt}[G_2(t)], \ldots .$$

Similarly,

$$e^{-A(t)\lambda^2}\widetilde{h}_M(\lambda,\Omega(\lambda,\cdot),t)=\frac{1}{\lambda^2}e^{-A(t)\lambda^2}\int_{\tau=0}^{t}\frac{d}{d\tau}[e^{A(\tau)\lambda^2}]\frac{e^{-B(\tau)}h_M(\lambda,\tau)}{\alpha(\tau)}d\tau$$

$$=\frac{1}{\lambda^2}e^{-A(t)\lambda^2}\int_{\tau=0}^{t}\frac{d}{d\tau}[e^{A(\tau)\lambda^2}]H_{M,0}(\lambda,\tau)d\tau$$

$$=\frac{1}{\lambda^2}H_{M,0}(\lambda,t)-\frac{1}{\lambda^2}e^{-A(t)\lambda^2}H_{M,0}(\lambda,0)-\frac{1}{\lambda^2}e^{-A(t)\lambda^2}\int_{\tau=0}^{t}e^{A(\tau)\lambda^2}\frac{d}{d\tau}[H_{M,0}(\lambda,\tau)]d\tau\,,$$

where

$$H_{M,0}(\lambda,\tau)=\frac{e^{-B(\tau)}h_M(\lambda,\tau)}{\alpha(\tau)}.$$

More generally, for fixed $M$ ,

$$\text{(3.26)}\quad e^{-A(t)\lambda^2}\widetilde{h}_M(\lambda,\Omega(\lambda,\cdot),t)=\sum_{s=1}^{K}(-1)^{s-1}\left[\frac{1}{\lambda^{2s}}H_{M,s-1}(\lambda,t)-\frac{1}{\lambda^{2s}}e^{-A(t)\lambda^2}H_{M,s}(\lambda,0)\right]$$

$$+(-1)^K\frac{1}{\lambda^{2K}}e^{-A(t)\lambda^2}\int_{\tau=0}^{t}e^{A(\tau)\lambda^2}\frac{d}{d\tau}[H_{M,K}(\lambda,\tau)]d\tau\quad(K\in\mathbb{N}),$$

where

$$H_{M,1}(t)=\frac{1}{\alpha(t)}\frac{d}{dt}[H_{M,0}(t)]\,,\ H_{M,2}(t)=\frac{1}{\alpha(t)}\frac{d}{dt}[H_{M,1}(t)]\,,\ H_{M,3}(t)=\frac{1}{\alpha(t)}\frac{d}{dt}[H_{M,2}(t)]\,,\,\ldots\,.$$

***Step 4*** We write $\mathcal{G}^t(g_0)(x,t)=\int_{\Gamma_0}+\int_{\Gamma_1}$ and, in view of (3.25), we have

$$\text{(3.27)}\quad \int_{\Gamma_1}\cdots d\lambda=-2ie^{B(t)}\sum_{l=1}^{L}(-1)^{l-1}\int_{\Gamma_1}\left[G_{l-1}(t)\frac{e^{i\lambda x}}{\lambda^{2l-1}}-G_{l-1}(0)\frac{e^{i\lambda x}}{\lambda^{2l-1}}e^{-A(t)\lambda^2}\right]d\lambda$$

$$-2i(-1)^Le^{B(t)}\left(\int_{\{\lambda\in\mathbb{R}:|\lambda\geq1\}|}-\int_{\gamma_1+\gamma_2}\right)\left[\frac{e^{i\lambda x}}{\lambda^{2L-1}}e^{-A(t)\lambda^2}\int_{\tau=0}^{t}e^{A(\tau)\lambda^2}\frac{d}{d\tau}[G_{L-1}(t)]d\tau\right]d\lambda\,.$$

Also,

$$\text{(3.28)}\quad \int_{\Gamma_1}G_{l-1}(t)\frac{e^{i\lambda x}}{\lambda^{2l-1}}d\lambda=\int_{\{|\lambda|=1\ \&\ \frac{\pi}{4}\leq\arg\lambda\leq\frac{3\pi}{4}\}}G_{l-1}(t)\frac{e^{i\lambda x}}{\lambda^{2l-1}}d\lambda$$

and

$$\text{(3.29)}\quad \int_{\Gamma_1}G_{l-1}(0)\frac{e^{i\lambda x}}{\lambda^{2l-1}}e^{-A(t)\lambda^2}d\lambda=\left(\int_{\{\lambda\in\mathbb{R}:|\lambda|\geq1\}}-\int_{\gamma_1+\gamma_2}\right)G_{l-1}(0)\frac{e^{i\lambda x}}{\lambda^{2l-1}}e^{-A(t)\lambda^2}d\lambda\,.$$

Substituting (3.28) and (3.29) in (3.27), we can easily prove the 2nd assertion for the integral $\mathcal{G}^t(g_0)$. The point here is that the integral in the RHS of (3.28) is taken over a finite contour, while the integral in the RHS of (3.29), which is taken over $\{\lambda\in\mathbb{R}:|\lambda|\geq1\}$, contains the factor $e^{-A(t)\lambda^2}$. These facts allow the interchange of differentiation (with respect to $x$ and $t$) and the limit (as $x\to0^+$) with the integrations (in the RHS (3.27)), as long as $t$ is kept away from $0$. (We remind the analogous argument in Step 2.)

***Step 5*** In order to prove the 2nd assertion for the integral $\mathcal{F}^{+}(f)$, we substitute (3.26) in (3.4) and observe that

$$(3.30)\quad \int_{\Gamma_1} \frac{e^{i\lambda x}}{\lambda^{2s-1}} H_{M,s-1}(\lambda,t)d\lambda = \int_{\{|\lambda|=1\ \&\ \frac{\pi}{4}\leq \arg\lambda \leq \frac{3\pi}{4}\}} \frac{e^{i\lambda x}}{\lambda^{2s-1}} H_{M,s-1}(\lambda,t)d\lambda$$

and

$$(3.31)\quad \left(\int_{\Gamma_1} + \int_{\gamma_1+\gamma_2}\right) \frac{e^{i\lambda x}}{\lambda^{2s-1}} e^{-A(t)\lambda^2} H_{M,s}(\lambda,0)d\lambda = \int_{\{\lambda\in\mathbb{R}:|\lambda\geq 1\}|} \frac{e^{i\lambda x}}{\lambda^{2s-1}} e^{-A(t)\lambda^2} H_{M,s}(\lambda,0)d\lambda .$$

The proof of the 2nd assertion for the integral $\mathcal{F}^{-}(f)$ is similar.
This completes the proof of the 2nd assertion.

***Step 6*** Having proved the 1st and 2nd assertions, the remaing conclusions of the theorem follow easily, once we have the results of Theorem 1.

**Proof of Theorem 3** Having at our disposal the expansions of the kernels of formula (2.9), as these are described by (3.1), (3.4), (3.25) and (3.26), it is easy to see that rather straightforward computations lead to the conclusions of the theorem.

**Proof of Theorem 4** ***Step 1*** *If the smooth function* $H(\lambda,t)$ *satisfies*

$$\sup_{t\leq T}\left|\frac{\partial^m H(\lambda,t)}{\partial t^m}\right| = \mathrm{O}(\lambda^N),\ \textit{as}\ \lambda\to\infty\ \textit{with}\ \lambda\in\Gamma_1,$$

*for some* $m$, $T>0$ *and* $N\in\mathbb{N}$, *then*

$$\lim_{x\to+\infty}\left\{x^{\ell}\frac{\partial^{n+m}}{\partial x^n\partial t^m}\left[\int_{\Gamma_1} e^{i\lambda x}H(\lambda,t)d\lambda\right]\right\}=0,$$

*for nonnegative integers* $n$ *and* $\ell$, *uniformly for* $0<t\leq T$.

*Proof* It suffices to notice that, for $x\geq 1$,

$$\left|x^{\ell}\frac{\partial^{n+m}}{\partial x^n\partial t^m}\left[\int_{\Gamma_1} e^{i\lambda x}H(\lambda,t)d\lambda\right]\right| = \left|x^{\ell}\int_{\Gamma_1}(i\lambda)^n e^{i\lambda x}\frac{\partial^m H(\lambda,t)}{\partial t^m}d\lambda\right|$$

$$\leq x^{\ell}e^{-[x\sin(\pi/4)]/2}\int_{\Gamma_1}|\lambda|^n e^{-[|\lambda|\sin(\pi/4)]/2}\left|\frac{\partial^m H(\lambda,t)}{\partial t^m}\right| d|\lambda| .$$

***Step 2*** *The functions defined by integrals taken on* $\Gamma$

Writing the integral $\mathcal{G}^t(g_0)=\int_{\Gamma}\cdots d\lambda$, which defines $\mathcal{G}^t(g_0)$, as

$$\mathcal{G}^t(g_0)=\int_{\Gamma_1}+\int_{\Gamma_0},$$

and applying the result of *Step1* with

$$H(\lambda,t) := e^{-\Omega(\lambda,t)}\int_{\tau=0}^{t} e^{\Omega(\lambda,\tau)}\alpha(\tau)g_0(\tau)d\tau ,$$

we see, firstly, that the part of $\mathcal{G}^t(g_0)(x,t)$, which comes from the integral over $\Gamma_1$, satisfies the conclusion of the theorem. To prove the corresponding assertion for the part of $\mathcal{G}^t(g_0)(x,t)$, which comes from the integral over $\Gamma_0$, we can integrate by parts several times – depending on $\ell$. For example, in the case $\ell = 1$, it suffices to write

$$x\frac{\partial^{n+m}}{\partial x^n \partial t^m}\left[\int_{\Gamma_0} i\lambda\left\{e^{i\lambda x-\Omega(\lambda,t)}\int_{\tau=0}^{t} e^{\Omega(\lambda,\tau)}\alpha(\tau)g_0(\tau)d\tau\right\}d\lambda\right]$$

$$=\frac{1}{i}\int_{\Gamma_0}(i\lambda)^{n+1}\left[\frac{d}{d\lambda}(e^{i\lambda x})\right]\left(\frac{\partial^m}{\partial t^m}\left\{e^{-\Omega(\lambda,t)}\int_{\tau=0}^{t} e^{\Omega(\lambda,\tau)}\alpha(\tau)g_0(\tau)d\tau\right\}\right)d\lambda$$

$$=(i\lambda)^{n+1}e^{i\lambda x}\left(\frac{\partial^m}{\partial t^m}\left\{e^{-\Omega(\lambda,t)}\int_{\tau=0}^{t} e^{\Omega(\lambda,\tau)}\alpha(\tau)g_0(\tau)d\tau\right\}\right)\Bigg|_{\lambda=e^{i3\pi/4}}^{\lambda=e^{i\pi/4}}$$

$$-\frac{1}{i}\int_{\Gamma_0} e^{i\lambda x}\frac{d}{d\lambda}\left((i\lambda)^{n+1}\frac{\partial^m}{\partial t^m}\left\{e^{-\Omega(\lambda,t)}\int_{\tau=0}^{t} e^{\Omega(\lambda,\tau)}\alpha(\tau)g_0(\tau)d\tau\right\}\right)d\lambda$$

and observe that the integrand of the last integral tends to zero, as $x \to +\infty$, for $\lambda \in \Gamma_0 - \{0\}$.

The integrals $\mathcal{U}^-(u_0)$ and $\mathcal{F}^-(f)$ can be treated similarly.

***Step 3*** *The function defined by the integral* $\mathcal{U}^+(u_0)$.

In this case we have to show that

$$\lim_{x\to+\infty}\left[x^{\ell}\frac{\partial^{n+m}\mathcal{U}^+(u_0)(x,t)}{\partial x^n \partial t^m}\right]=0. \tag{3.32}$$

If $\ell = 0$, (3.32) follows from the representstion (3.3): For the first and the third integral in the RHS of (3.3), we use the Riemann-Lebesgue lemma, while for the second integral in (3.3), we work as in *Step 1*.

If $\ell = 1$, we use again (3.3) (with sufficiently large $M$) and we integrate by parts and we obtain

$$x\frac{\partial^{n+m}\mathcal{U}^+(u_0)(x,t)}{\partial x^n \partial t^m}=-\frac{1}{i}\int_{-1}^{1} e^{i\lambda x}\frac{d}{d\lambda}\left\{(i\lambda)^n\frac{\partial^m[e^{-\Omega(\lambda,t)}]}{\partial t^m}\hat{u}_0(\lambda)\right\}d\lambda \tag{3.33}$$

$$-\frac{1}{i}\sum_{k=1}^{M}u_0^{(k-1)}(0)\left(\int_{\Gamma_1}+\int_{\gamma_1+\gamma_2}\right)e^{i\lambda x}\frac{d}{d\lambda}\left\{(i\lambda)^n\frac{\partial^m[e^{-\Omega(\lambda,t)}]}{\partial t^m}\frac{1}{(i\lambda)^k}\right\}d\lambda$$

$$-\frac{1}{i}\left(\int_{-\infty}^{-1}+\int_{1}^{\infty}\right)e^{i\lambda x}\frac{d}{d\lambda}\left\{(i\lambda)^n\frac{\partial^m[e^{-\Omega(\lambda,t)}]}{\partial t^m}(u_0^{(M)}\hat{)}(\lambda)\frac{1}{(i\lambda)^M}\right\}d\lambda,$$

since the "intermediate evaluations" (boundary terms) cancel each other out. Then for the first and the third integral in the RHS of (3.33), we use the Riemann-Lebesgue lemma, while for the second integral in (3.33), we work as in *Step 1*.

Further integration by parts, proves (3.32) for the general case: $\ell \ge 2$.

***Step 4*** *The function defined by the integral* $\mathcal{F}^+(f)$.

In this case we have to show that

(3.34) $$\lim_{x\to+\infty}\left[x^{\ell}\,\frac{\partial^{n+m}\mathcal{F}^{+}(f)(x,t)}{\partial x^{n}\partial t^{m}}\right]=0\,.$$

Using (3.5), with sufficiently large $M$, the proof of (3.34) is similar to the proof of (3.32).

This completes the proof of the theorem.

***Remark*** Regarding the proof of Theorem 4, we insisted on a detailed analysis of the integrals, since in cases of other equations the corresponding UTM solution need not be rapidly decreasing as $x\to+\infty$, despite the relevant assumption on the datum $u_0$. See for example [19] and [32].

### 4. The Barenblatt equation: Derivation of the solution and statement of the results

The formal construction of the UTM-type solution formula in this case has been carried out in [26], however we present it here again for the convenience of readers. We will then verify it and analyse it. Working with $\lambda\in\mathbb{C}$, $\operatorname{Im}\lambda\le 0$, and taking Fourier transforms, the differential equation in (1.2) gives

(4.1) $$\frac{\partial}{\partial t}[\hat{u}(\lambda,t)]=\alpha(t)\frac{\partial}{\partial t}\left[\int_0^\infty\frac{\partial^2 u(y,t)}{\partial y^2}e^{-i\lambda y}dy\right]+\beta(t)\int_0^\infty\frac{\partial^2 u(y,t)}{\partial y^2}e^{-i\lambda y}dy$$
$$+\gamma(t)\int_0^\infty u(y,t)e^{-i\lambda y}dy+\int_0^\infty f(y,t)e^{-i\lambda y}dy\,.$$

Also,

$$\int_0^\infty u_x(x,t)e^{-i\lambda x}dx=-g_0(t)+i\lambda\hat{u}(\lambda,t)\ \text{ and }\ \int_0^\infty u_{xx}(x,t)e^{-i\lambda x}dx=-g_1(t)-i\lambda g_0(t)+(i\lambda)^2\hat{u}(\lambda,t)\,,$$

where we have set $g_1(t)=u_x(0,t)$, so, integrating by parts in (4.1), we obtain:

$$[1+\alpha(t)\lambda^2]\frac{\partial}{\partial t}[\hat{u}(\lambda,t)]+[\beta(t)\lambda^2-\gamma(t)]\hat{u}(\lambda,t)$$
$$=-\alpha(t)[g_1'(t)+i\lambda g_0'(t)]-\beta(t)[g_1(t)+i\lambda g_0(t)]+\hat{f}(\lambda,t)\,,$$

where we used also $\frac{\partial}{\partial t}[u_x(0,t)]=g_1'(t)$.

Thus,

(4.2) $$\frac{\partial}{\partial t}[\hat{u}(\lambda,t)]+\frac{\beta(t)\lambda^2-\gamma(t)}{1+\alpha(t)\lambda^2}\hat{u}(\lambda,t)$$
$$=-\frac{1}{1+\alpha(t)\lambda^2}\{\alpha(t)[g_1'(t)+i\lambda g_0'(t)]+\beta(t)[g_1(t)+i\lambda g_0(t)]\}+\frac{\hat{f}(\lambda,t)}{1+\alpha(t)\lambda^2}\,.$$

Setting

$$\phi(\lambda,t)=\frac{\beta(t)\lambda^2-\gamma(t)}{1+\alpha(t)\lambda^2}\ \text{ and }\ \Phi(\lambda,t)=\int_{\tau=0}^{t}\phi(\lambda,\tau)d\tau=\int_{\tau=0}^{t}\frac{\beta(\tau)\lambda^2-\gamma(\tau)}{1+\alpha(\tau)\lambda^2}d\tau\,,$$

we see – and this is crucial for what follows – that

(4.3) $$\phi(-\lambda,t)=\phi(\lambda,t)\,,\ \Phi(-\lambda,t)=\Phi(\lambda,t)\,,$$

and

(4.4) $$\lim_{\substack{\lambda\to\infty \\ \lambda\in\mathbb{C}-i\mathbb{R}}} \Phi(\lambda,t) = \int_{\tau=0}^{t} \frac{\beta(\tau)}{\alpha(\tau)} d\tau$$, uniformly for $t$ in compact subsets of $[0,+\infty)$.

Assuming

(4.5) $$c_0 := \inf\{1/\sqrt{\alpha(t)} : t \ge 0\} > 0\} < +\infty ,$$

it follows from (4.4) that, for every $T > 0$ and for every open neighborhood $\mathcal{A}$ of the half-line $\ell_0 := [c_0 i, +\infty i) := \{\eta i : \eta \ge c_0\}$, such that $\inf_{\mathrm{E}} lengh(\mathcal{A} \cap \mathrm{E}) > 0$ (where $\mathrm{E}$ runs through the straight lines perpendicular to the $\eta - axis$ at the points of $\ell_0$),

(4.6) $$\sup\{ \left| e^{\pm\Phi(\lambda,t)} \right| : 0 \le t \le T,\ \lambda \in \mathbb{C} - \mathcal{A} \ \ and \ \ \mathrm{Im}\,\lambda \ge 0\} < +\infty .$$

We also point out that, as a straightforward computation shows,

(4.7) $$\mathfrak{P}[e^{i\lambda x - \Phi(\lambda,t)}] = 0 ,$$

where $\mathfrak{P}$ is the differential operator

$$\mathfrak{P} = \mathfrak{P}_{x,t} = \frac{\partial}{\partial t} - \alpha(t)\frac{\partial^3}{\partial t \partial x^2} - \beta(t)\frac{\partial^2}{\partial x^2} - \gamma(t) ,$$

acting on functions defined for $(x,t) \in Q = \mathbb{R}^+ \times \mathbb{R}^+$. Thus, with this notation, the differential equation in (1.2) can be written as $\mathfrak{P}u = f$.

Now, since

$$\frac{\partial}{\partial t}[e^{\Phi(\lambda,t)}\hat{u}(\lambda,t)] = e^{\Phi(\lambda,t)} \frac{\partial}{\partial t}[\hat{u}(\lambda,t)] + e^{\Phi(\lambda,t)} \frac{\beta(t)\lambda^2 - \gamma(t)}{1+\alpha(t)\lambda^2}\hat{u}(\lambda,t) ,$$

we may write (4.2) as follows:

$$\frac{\partial}{\partial t}[e^{\Phi(\lambda,t)}\hat{u}(\lambda,t)] = -\frac{e^{\Phi(\lambda,t)}}{1+\alpha(t)\lambda^2}\{\alpha(t)[g_1'(t) + i\lambda g_0'(t)] + \beta(t)[g_1(t) + i\lambda g_0(t)]\} + \frac{e^{\Phi(\lambda,t)}}{1+\alpha(t)\lambda^2}\hat{f}(\lambda,t) .$$

Integrating the above equation, we obtain

$$e^{\Phi(\lambda,t)}\hat{u}(\lambda,t) = \hat{u}_0(\lambda) - \int_{\tau=0}^{t} \frac{e^{\Phi(\lambda,\tau)}}{1+\alpha(\tau)\lambda^2}\left\{\alpha(\tau)[g_1'(\tau) + i\lambda g_0'(\tau)] + \beta(\tau)[g_1(\tau) + i\lambda g_0(\tau)] - \hat{f}(\lambda,\tau)\right\} d\tau ,$$

or equivalently,

(4.8) $$\hat{u}(\lambda,t) = e^{-\Phi(\lambda,t)}\hat{u}_0(\lambda) - e^{-\Phi(\lambda,t)} \int_{\tau=0}^{t} \frac{e^{\Phi(\lambda,\tau)}}{1+\alpha(\tau)\lambda^2}[\alpha(\tau)g_1'(\tau) + \beta(\tau)g_1(\tau)]\, d\tau$$
$$- e^{-\Omega(\lambda,t)} i\lambda \int_{\tau=0}^{t} \frac{e^{\Phi(\lambda,\tau)}}{1+\alpha(\tau)\lambda^2}[\alpha(\tau)g_0'(\tau) + \beta(\tau)g_0(\tau)]\, d\tau$$
$$+ e^{-\Phi(\lambda,t)} \int_{\tau=0}^{t} \frac{e^{\Phi(\lambda,\tau)}}{1+\alpha(\tau)\lambda^2}\hat{f}(\lambda,\tau)\, d\tau ,$$ for $\lambda \in \mathbb{C} - [-ic_0, -\infty i)$ with $\mathrm{Im}\,\lambda \le 0$.

In view of (4.3), setting $-\lambda$ in place of $\lambda$, in (4.8), we find

$$(4.9)\quad \hat{u}(-\lambda,t) = e^{-\Phi(\lambda,t)}\hat{u}_0(-\lambda) - e^{-\Phi(\lambda,t)} \int_{\tau=0}^{t} \frac{e^{\Phi(\lambda,\tau)}}{1+\alpha(\tau)\lambda^2}[\alpha(\tau)g_1'(\tau)+\beta(\tau)g_1(\tau)]\,d\tau$$

$$- e^{-\Phi(\lambda,t)} \int_{\tau=0}^{t} \frac{e^{\Phi(\lambda,\tau)}}{1+\alpha(\tau)\lambda^2}\left\{-i\lambda[\alpha(\tau)g_0'(\tau)+\beta(\tau)g_0(\tau)] - \hat{f}(-\lambda,\tau)\right\}d\tau, \text{ for } \lambda\in\mathbb{C} \text{ with } \operatorname{Im}\lambda \geq 0.$$

Multiplying (4.8) by $e^{i\lambda x}$ and integrating with respect to $\lambda$, we obtain

$$(4.10)\quad 2\pi u(x,t) = \int_{-\infty}^{\infty} e^{i\lambda x}\hat{u}(\lambda,t)\,d\lambda$$

$$= \int_{-\infty}^{\infty} e^{i\lambda x-\Phi(\lambda,t)}\hat{u}_0(\lambda)\,d\lambda - \int_{-\infty}^{\infty}\left\{ e^{i\lambda x-\Phi(\lambda,t)} \int_{\tau=0}^{t} \frac{e^{\Phi(\lambda,\tau)}}{1+\alpha(\tau)\lambda^2}[\alpha(\tau)g_1'(\tau)+\beta(\tau)g_1(\tau)]\,d\tau\right\}d\lambda$$

$$- \int_{-\infty}^{\infty}\left\{ i\lambda e^{i\lambda x-\Phi(\lambda,t)} \int_{\tau=0}^{t} \frac{e^{\Phi(\lambda,\tau)}}{1+\alpha(\tau)\lambda^2}[\alpha(\tau)g_0'(\tau)+\beta(\tau)g_0(\tau)]\,d\tau\right\}d\lambda$$

$$+ \int_{-\infty}^{\infty}\left\{ e^{i\lambda x-\Phi(\lambda,t)} \int_{\tau=0}^{t} \frac{e^{\Phi(\lambda,\tau)}}{1+\alpha(\tau)\lambda^2}\hat{f}(\lambda,\tau)\,d\tau\right\}d\lambda.$$

Using (4.6) and deforming the contours of the $d\lambda$-integrals which involve the quantities

$$\alpha(\tau)g_1'(\tau)+\beta(\tau)g_1(\tau),\ \alpha(\tau)g_0'(\tau)+\beta(\tau)g_0(\tau),$$

we write (4.10) in the following way:

$$(4.11)\quad 2\pi u(x,t) = \int_{-\infty}^{\infty} e^{i\lambda x}\hat{u}(\lambda,t)\,d\lambda$$

$$= \int_{-\infty}^{\infty} e^{i\lambda x-\Phi(\lambda,t)}\hat{u}_0(\lambda)\,d\lambda - \int_{\mathcal{L}}\left\{ e^{i\lambda x-\Phi(\lambda,t)} \int_{\tau=0}^{t} \frac{e^{\Phi(\lambda,\tau)}}{1+\alpha(\tau)\lambda^2}[\alpha(\tau)g_1'(\tau)+\beta(\tau)g_1(\tau)]\,d\tau\right\}d\lambda$$

$$- \int_{\mathcal{L}}\left\{ i\lambda e^{i\lambda x-\Phi(\lambda,t)} \int_{\tau=0}^{t} \frac{e^{\Phi(\lambda,\tau)}}{1+\alpha(\tau)\lambda^2}[\alpha(\tau)g_0'(\tau)+\beta(\tau)g_0(\tau)]\,d\tau\right\}d\lambda$$

$$+ \int_{-\infty}^{\infty}\left\{ e^{i\lambda x-\Phi(\lambda,t)} \int_{\tau=0}^{t} \frac{e^{\Phi(\lambda,\tau)}}{1+\alpha(\tau)\lambda^2}\hat{f}(\lambda,\tau)\,d\tau\right\}d\lambda,$$

where $\mathcal{L}$ is a contour as the one depicted in fig.3, below. More precisely, $\mathcal{L}$ consists of two half-lines, parallel to the imaginary axis, and a half-circle-like curve, joining the half-lines, and oriented counterclockwise, as in fig. 3. "Minor" deformations of the half-lines will not change the values of the integrals in (4.11), which are taken over such a deformed $\mathcal{L}$. For, example $\mathcal{L}$ can be replaced in (4.11) by $\partial\mathcal{A}$, where $\mathcal{A}$ is as in (4.6), provided that the contour $\partial\mathcal{A}$ is sufficiently smooth.

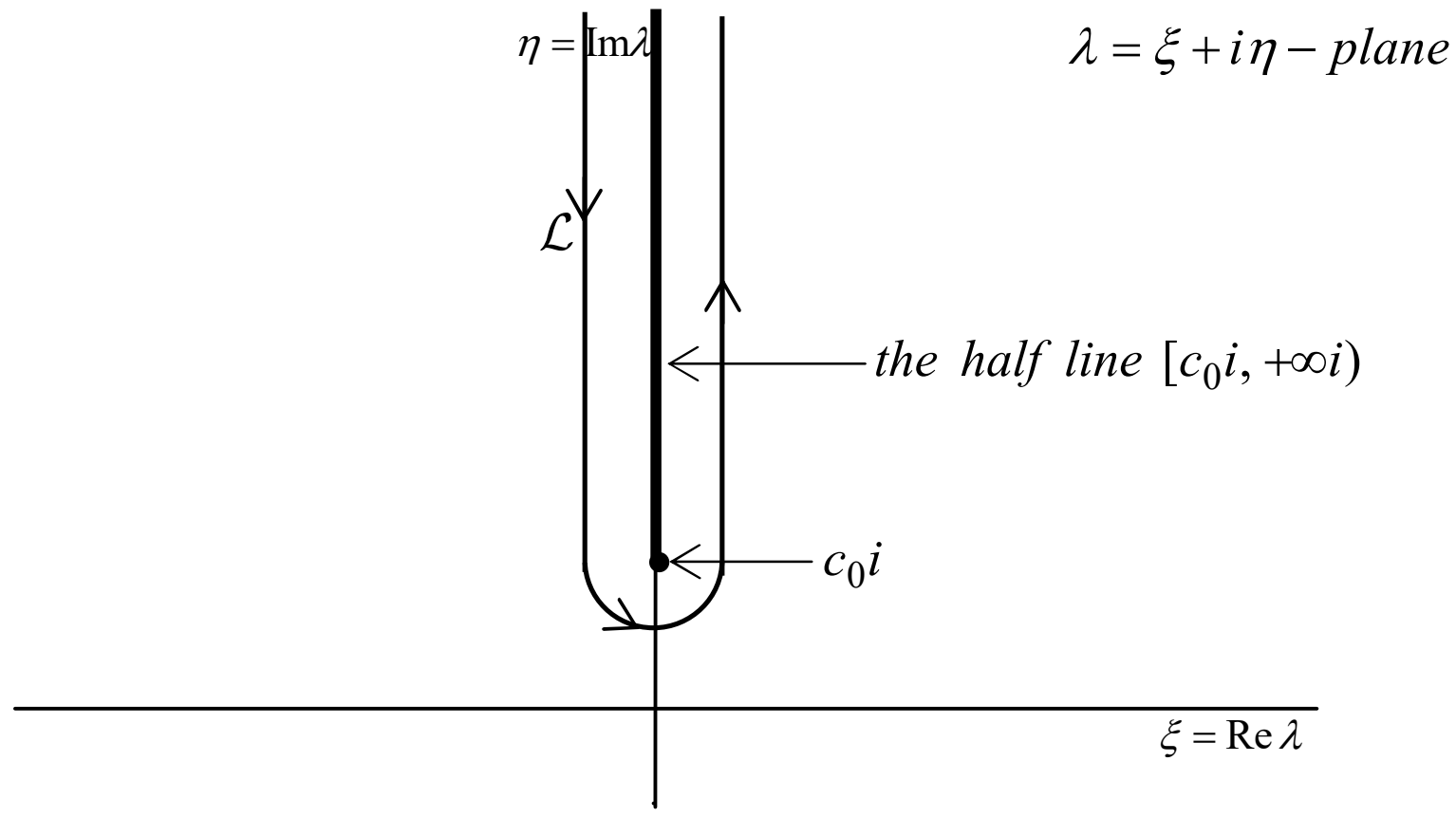

**Fig. 3** A choice of the contour $\mathcal{L}$.

On the other hand, (4.9) gives

$$(4.12)\quad 0 = \int_{\mathcal{L}} e^{i\lambda x} \hat{u}(-\lambda,t)d\lambda$$

$$= \int_{\mathcal{L}} e^{i\lambda x-\Phi(\lambda,t)} \hat{u}_0(-\lambda)d\lambda - \int_{\mathcal{L}} \left\{ e^{i\lambda x-\Phi(\lambda,t)} \int_{\tau=0}^{t} \frac{e^{\Phi(\lambda,\tau)}}{1+\alpha(\tau)\lambda^2} [\alpha(\tau)g_1'(\tau) + \beta(\tau)g_1(\tau)]\, d\tau \right\} d\lambda$$

$$+ \int_{\mathcal{L}} \left\{ i\lambda e^{i\lambda x-\Phi(\lambda,t)} \int_{\tau=0}^{t} \frac{e^{\Phi(\lambda,\tau)}}{1+\alpha(\tau)\lambda^2} [\alpha(\tau)g_0'(\tau) + \beta(\tau)g_0(\tau)]\, d\tau \right\} d\lambda$$

$$+ \int_{\mathcal{L}} \left\{ e^{i\lambda x-\Phi(\lambda,t)} \int_{\tau=0}^{t} \frac{e^{\Phi(\lambda,\tau)}}{1+\alpha(\tau)\lambda^2} \hat{f}(-\lambda,\tau)\, d\tau \right\} d\lambda .$$

Subtracting (4.12) from (4.11), we obtain the solution to *Problem 2*: For $x > 0$ and $t > 0$,

$$(4.13)\quad 2\pi u(x,t) = \int_{-\infty}^{\infty} e^{i\lambda x-\Phi(\lambda,t)} \hat{u}_0(\lambda)d\lambda - \int_{\mathcal{L}} e^{i\lambda x-\Phi(\lambda,t)} \hat{u}_0(-\lambda)d\lambda$$

$$-2\int_{\mathcal{L}} \left\{ i\lambda e^{i\lambda x-\Phi(\lambda,t)} \int_{\tau=0}^{t} \frac{e^{\Phi(\lambda,\tau)}}{1+\alpha(\tau)\lambda^2} [\alpha(\tau)g_0'(\tau) + \beta(\tau)g_0(\tau)]\, d\tau \right\} d\lambda$$

$$+ \int_{-\infty}^{\infty} \left\{ e^{i\lambda x-\Phi(\lambda,t)} \int_{\tau=0}^{t} \frac{e^{\Phi(\lambda,\tau)}}{1+\alpha(\tau)\lambda^2} \hat{f}(\lambda,\tau)\, d\tau \right\} d\lambda - \int_{\mathcal{L}} \left\{ e^{i\lambda x-\Phi(\lambda,t)} \int_{\tau=0}^{t} \frac{e^{\Phi(\lambda,\tau)}}{1+\alpha(\tau)\lambda^2} \hat{f}(-\lambda,\tau)\, d\tau \right\} d\lambda .$$

***Remarks*** **(1)** We point out the following: If $T > 0$ is fixed and $t$ is restricted to the interval $0 < t < T$, then the contour $\mathcal{L}$ in the above formula can be replaced by a finite contour $\mathcal{C}_T$ of the following form: $\mathcal{C}_T$ can be a simple closed curve surrounding the line segment $[c_0 i, c_T i] = \{\eta i : c_0 \le \eta \le c_T\}$, where $c_T := \sup\{1/\sqrt{\alpha(t)} : 0 \le t \le T\}$, and oriented counterclockwise, as depicted in fig. 4, below.

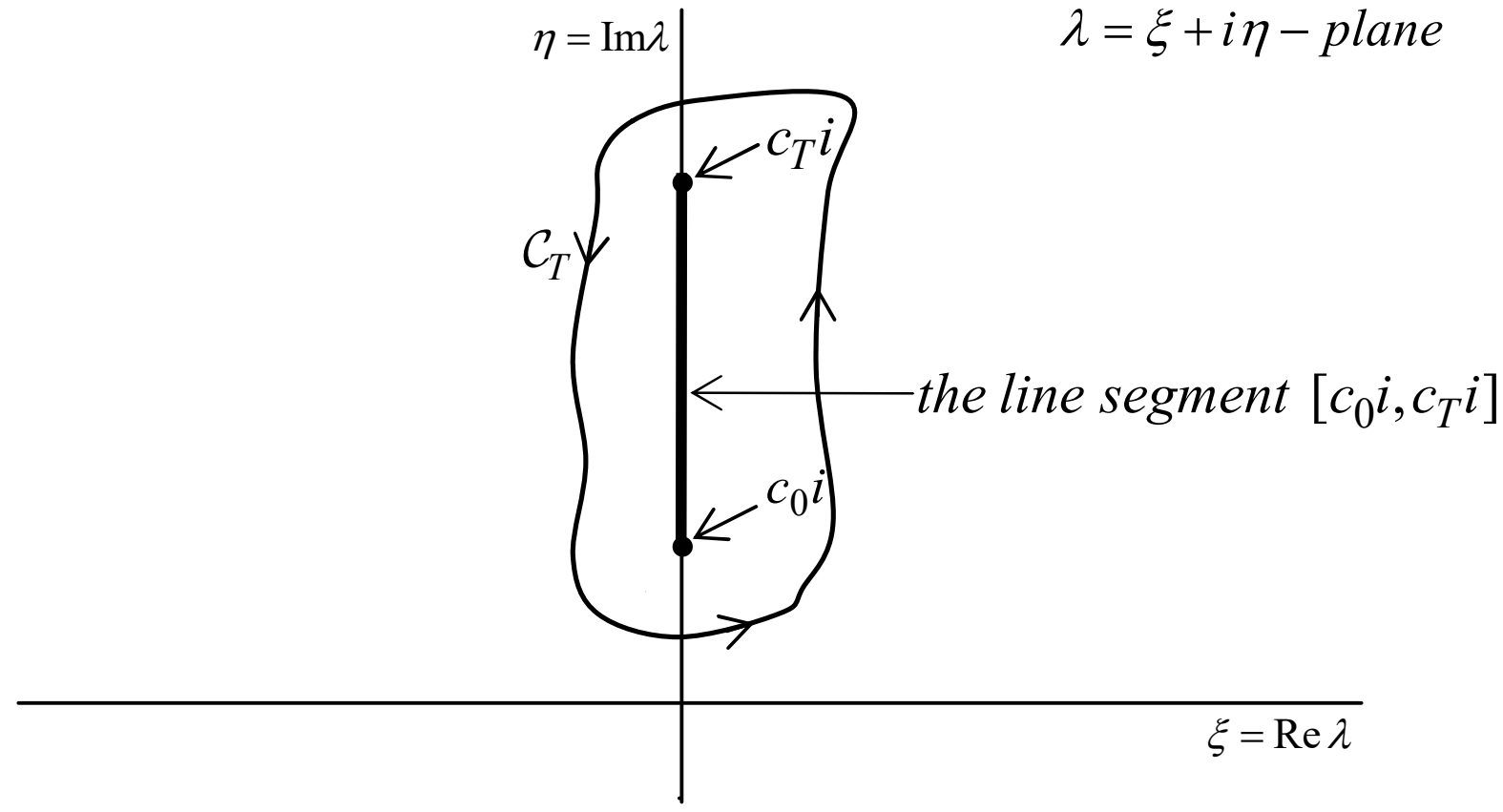

**Fig. 4** A choice of the contour $\mathcal{C}_T$.

This will not alter the values of the corresponding integrals. Thus, in many parts of what follows a contour like $\mathcal{C}_T$ can and will be used in place of $\mathcal{L}$, provided that $t$ is restricted to the interval $0 < t < T$. This will be convenient and will simplify certain of our computations.

**(2)** In the case the function $\alpha(t)$ satisfies, in addition to (4.5), the condition

$$c_\infty := \sup\{1/\sqrt{\alpha(t)} : t \geq 0\} < +\infty ,$$

the contour $\mathcal{L}$, in formula (4.11) (and (4.13)), can be replaced by any finite closed contour $\mathcal{C}$, surrounding the line segment $[c_0 i, c_\infty i] = \{\eta i : \; c_0 \leq \eta \leq c_\infty\}$ and oriented counterclockwise, as depicted in fig.5, below. (This is essentially a special case of the previous remark.)

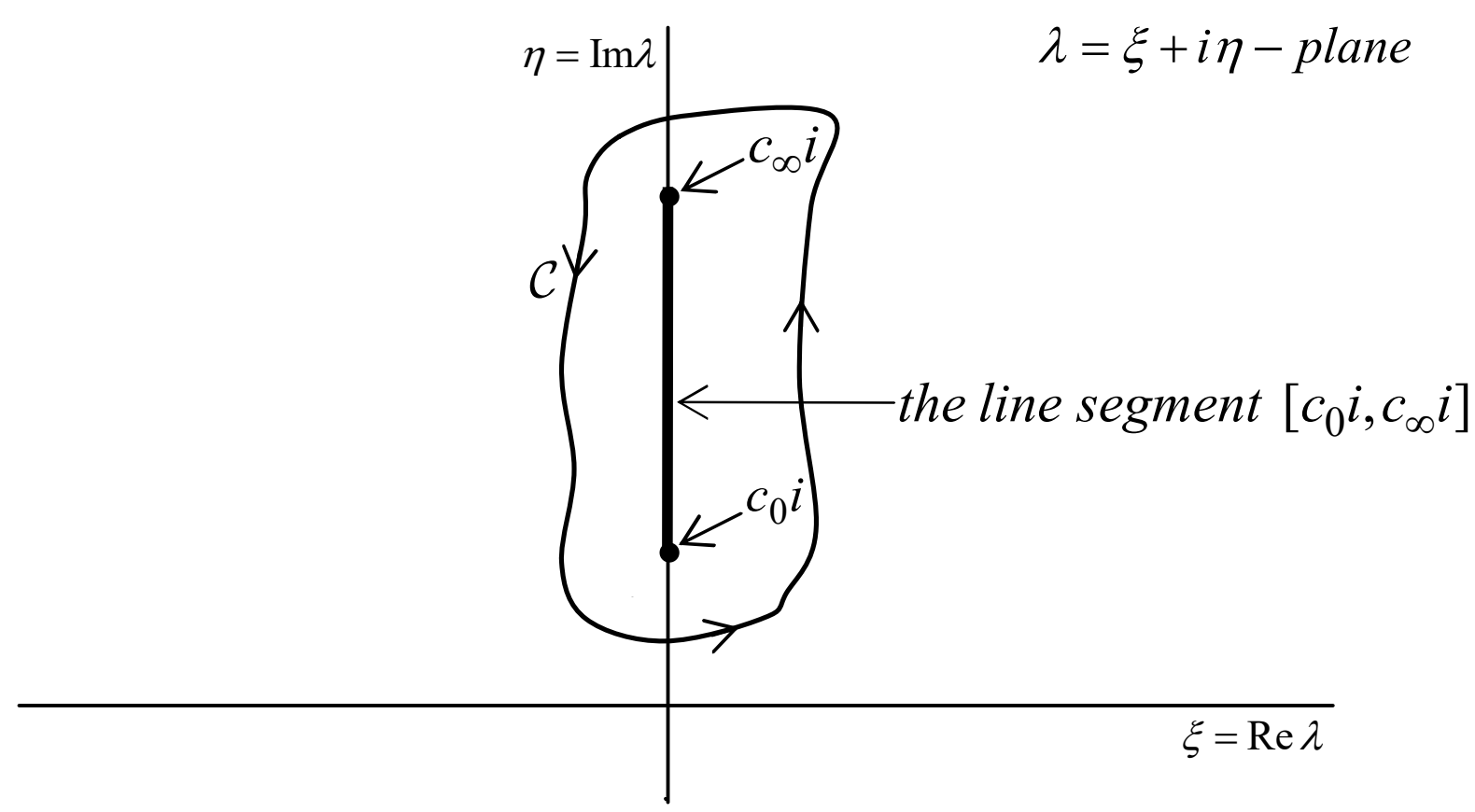

**Fig**. **5** A choice of the contour $\mathcal{C}$.

***Notation*** We will donote by $\mathfrak{U}^+(x,t)$, $\mathfrak{U}^-(x,t)$, $\mathcal{I}_{g_0,g_0'}(x,t)$, $\mathfrak{F}^+(x,t)$ and $\mathfrak{F}^-(x,t)$ the functions defined by the $d\lambda$-integrals in (4.13) which involve the terms $\hat{u}_0(\lambda)$, $\hat{u}_0(-\lambda)$, $g_0(\tau) \,\&\, g_0'(\tau)$, $\hat{f}(\lambda,\tau)$ and $\hat{f}(-\lambda,\tau)$ , respectively.

Thus, (4.13) can be written as follows:

(4.14) $$2\pi u(x,t) = \mathfrak{U}^+(x,t) - \mathfrak{U}^-(x,t) - 2\mathcal{I}_{g_0,g_0'}(x,t) + \mathfrak{F}^+(x,t) - \mathfrak{U}^-(x,t), \; (x,t) \in Q .$$

Also, setting

$$\mathfrak{G}(\lambda,t)=\int_{\tau=0}^{t}\frac{e^{\Phi(\lambda,\tau)}}{1+\alpha(\tau)\lambda^2}[\alpha(\tau)g_0'(\tau)+\beta(\tau)g_0(\tau)]\,d\tau\,,$$

we have

$$\mathcal{I}_{g_0,g_0'}(x,t)=\int_{\mathcal{L}}\{i\lambda e^{i\lambda x-\Phi(\lambda,t)}\mathfrak{G}(\lambda,t)\}d\lambda\,. \tag{4.15}$$

*An equivalent form of* $\mathcal{I}_{g_0,g_0'}$ Integrating by parts, we have

$$\int_{\tau=0}^{t}\frac{e^{\Phi(\lambda,\tau)}}{1+\alpha(\tau)\lambda^2}\alpha(\tau)g_0'(\tau)\,d\tau=\frac{e^{\Phi(\lambda,\tau)}}{1+\alpha(\tau)\lambda^2}\alpha(\tau)g_0(\tau)\Bigg|_{\tau=0}^{\tau=t}-\int_{\tau=0}^{t}\frac{d}{d\tau}\left[\frac{e^{\Phi(\lambda,\tau)}}{1+\alpha(\tau)\lambda^2}\alpha(\tau)\right]g_0(\tau)\,d\tau$$

$$=\frac{e^{\Phi(\lambda,t)}}{1+\alpha(t)\lambda^2}\alpha(t)g_0(t)-\frac{\alpha(0)}{1+\alpha(0)\lambda^2}g_0(0)\;-\int_{\tau=0}^{t}\frac{e^{\Phi(\lambda,\tau)}}{1+\alpha(\tau)\lambda^2}\alpha'(\tau)g_0(\tau)\,d\tau$$

$$+\lambda^2\int_{\tau=0}^{t}\frac{e^{\Phi(\lambda,\tau)}}{[1+\alpha(\tau)\lambda^2]^2}\alpha(\tau)\alpha'(\tau)g_0(\tau)\,d\tau-\int_{\tau=0}^{t}\frac{e^{\Phi(\lambda,\tau)}}{1+\alpha(\tau)\lambda^2}\frac{\beta(\tau)\lambda^2-\gamma(\tau)}{1+\alpha(\tau)\lambda^2}\alpha(\tau)g_0(\tau)\,d\tau$$

$$=\frac{e^{\Phi(\lambda,t)}}{1+\alpha(t)\lambda^2}\alpha(t)g_0(t)-\frac{\alpha(0)}{1+\alpha(0)\lambda^2}g_0(0)\;-\int_{\tau=0}^{t}\frac{e^{\Phi(\lambda,\tau)}}{[1+\alpha(\tau)\lambda^2]^2}\alpha'(\tau)g_0(\tau)\,d\tau$$

$$-\int_{\tau=0}^{t}\frac{e^{\Phi(\lambda,\tau)}}{1+\alpha(\tau)\lambda^2}\frac{\beta(\tau)\lambda^2-\gamma(\tau)}{1+\alpha(\tau)\lambda^2}\alpha(\tau)g_0(\tau)\,d\tau\,.$$

Substituting in $\mathfrak{G}(\lambda,t)$, we see that (4.15) becomes:

$$\begin{aligned}\mathcal{I}_{g_0,g_0'}(x,t)&=\alpha(t)g_0(t)\int_{\mathcal{L}}\frac{i\lambda e^{i\lambda x}d\lambda}{1+\alpha(t)\lambda^2}-\alpha(0)g_0(0)\int_{\mathcal{L}}\frac{i\lambda e^{i\lambda x-\Phi(\lambda,t)}d\lambda}{1+\alpha(0)\lambda^2}\\&-\int_{\mathcal{L}}\left\{i\lambda e^{i\lambda x-\Phi(\lambda,t)}\int_{\tau=0}^{t}e^{\Phi(\lambda,\tau)}\frac{\beta(\tau)+\alpha(\tau)\gamma(\tau)}{[1+\alpha(\tau)\lambda^2]^2}g_0(\tau)\,d\tau\right\}d\lambda\\&-\int_{\mathcal{L}}\left\{i\lambda e^{i\lambda x-\Phi(\lambda,t)}\int_{\tau=0}^{t}\frac{e^{\Phi(\lambda,\tau)}}{[1+\alpha(\tau)\lambda^2]^2}\alpha'(\tau)g_0(\tau)\,d\tau\right\}d\lambda\\&+\int_{\mathcal{L}}\left\{i\lambda e^{i\lambda x-\Phi(\lambda,t)}\int_{\tau=0}^{t}\frac{e^{\Phi(\lambda,\tau)}}{[1+\alpha(\tau)\lambda^2]^2}\gamma(\tau)\alpha(\tau)g_0(\tau)\,d\tau\right\}d\lambda\,.\end{aligned} \tag{4.16}$$

Furthermore, computing the first integral in the RHS of (4.16), we find that, for fixed $x>0$ and $t>0$,

$$\int_{\mathcal{L}}\frac{i\lambda e^{i\lambda x}d\lambda}{1+\alpha(t)\lambda^2}=2\pi i\,Res\left[\frac{i\lambda e^{i\lambda x}}{1+\alpha(t)\lambda^2},\,\lambda=\frac{1}{i\sqrt{\alpha(t)}}\right]=-\frac{\pi}{\alpha(t)}e^{-x/\sqrt{\alpha(t)}}\,. \tag{4.17}$$

Thus, in view of (4.17), the “first” term in the quantity $-2\mathcal{I}_{g_0,g_0'}(x,t)$ (this quantity appears in the solution formula (4.14)) is equal to $2\pi e^{-x/\sqrt{\alpha(t)}}g_0(t)$. It is this term that will give the boundary condition $\lim_{x\to 0^+} u(x,t) = u(x,t)\big|_{x=0} = g_0(t)$. The presence of the other three terms (4.16) ensure that function $\mathcal{I}_{g_0,g_0'}(x,t)$ satisfies the homogeneous equation $\mathfrak{P}[\mathcal{I}_{g_0,g_0'}(x,t)] = 0$.

Now we state the main results for the solution of Problem 2. We recall our assumptions (1.3), (1.4), (1.9), (1.10) and (4.5).

**Theorem 5** *The function $u(x,t)$, defined by (4.13), is $C^\infty$ for $(x,t)\in Q = \mathbb{R}^+\times\mathbb{R}^+$, it satisfies the differential equation in (1.2), and*

*(4.18)*
$$\lim_{t\to 0^+} u(x,t) = u_0(x), \text{ uniformly for } x>0.$$

*Furthermore, if we assume, in addition, that $u_0(0) = g_0(0)$, then, for every $T>0$,*

*(4.19)*
$$\lim_{x\to 0^+} u(x,t) = g_0(t), \text{ uniformly for } 0<t\le T.$$

**Theorem 6** *Assuming (1.3), (1.4), (1.9), (1.10), and that $u_0(0) = g_0(0)$, the function $u(x,t)$, defined by (4.13), for $(x,t)\in Q = \mathbb{R}^+\times\mathbb{R}^+$, satisfies the following:*

*1st The limits*

$$G_{n,m}(t_0) := \lim_{\substack{(x,t)\to(0,\,t_0)\\(x,t)\in Q}} \frac{\partial^{n+m}u(x,t)}{\partial x^n \partial t^m},\ n,m = 0,1,2,...,$$

*exist and are uniform for $0<t_0\le T$ (for every $T>0$). Moreover the functions $G_{n,m}(t_0)$ are $C^\infty$ for $t_0\in[0,\infty)$.*

*2nd The limits*

$$U_{n,m}(x_0) := \lim_{\substack{(x,t)\to(x_0,\,t)\\(x,t)\in Q}} \frac{\partial^{n+m}u(x,t)}{\partial x^n \partial t^m},\ n = 0,1,2,...,$$

*exist and are uniform for $x_0>0$. Moreover the functions $U_{n,m}(x_0)$ are $C^\infty$ for $x_0\in[0,\infty)$.*

*3rd The function $u(x,t)$, extended to $\overline{Q}-\{(0,0)\}$ by setting*

$$u^*(x,0) := u_0(x) \text{ for } x>0 \text{ and } u^*(0,t) := g_0(t) \text{ for } t>0,$$

*is $C^\infty$ in $\overline{Q}-\{(0,0)\}$.*

*4th The limit $\lim_{t\to 0^+} \frac{\partial^n u(x,t)}{\partial x^n} = \frac{d^n u_0(x)}{dx^n}$, uniformly for $x>0$.*

*5th The limit $\lim_{x\to 0^+} \frac{\partial^m u(x,t)}{\partial t^m} = \frac{d^m g_0(t)}{dt^m}$, uniformly for $0<t\le T$.*

*6th The limits*

$$C_{n,m} := \lim_{\substack{(x,t)\to(0,0)\\(x,t)\in\overline{Q}-\{(0,0)\}}} \frac{\partial^{n+m}u^*(x,t)}{\partial x^n \partial t^m},\ n,m = 0,1,2,...,$$

*exist. Also $C_{n,0} = u_0^{(n)}(0)$ and $C_{0,m} = g_0^{(m)}(0)$.*

**Theorem 7** *With the assumptions of Theorem 6, the function* $u(x,t)$*, defined for* $(x,t)\in Q$ *by (4.13), satisfies, for* $T>0$*,*

$$\lim_{x\to+\infty}\left(x^{\ell}\,\frac{\partial^{n+m}u(x,t)}{\partial x^{n}\partial t^{m}}\right)=0,$$

*for nonnegative integers* $n$*,* $m$ *and* $\ell$*, uniformly for* $0<t\le T$*.*

## 5. Proofs of Theorems 5 - 7

**Proof of Theorem 5** Firstly, recalling formulas (1.7) and (1.8), and defining

$$\tilde{\hat{f}}(\lambda,\Phi(\lambda,\cdot),t)=\int_{\tau=0}^{t}\frac{e^{\Phi(\lambda,\tau)}}{1+\alpha(\tau)\lambda^2}\hat{f}(\lambda,\tau)\,d\tau,$$

we have

$$\tilde{\hat{f}}(\lambda,\Phi(\lambda,\cdot),t)=\int_0^t\frac{e^{\Phi(\lambda,\tau)}}{1+\alpha(\tau)\lambda^2}h_M(\lambda,\tau)d\tau+\int_0^t\frac{e^{\Phi(\lambda,\tau)}}{1+\alpha(\tau)\lambda^2}\delta_M(\lambda,\tau)d\tau \tag{5.1}$$

$$=:\tilde{h}_M(\lambda,\Phi(\lambda,\cdot),t)+\tilde{\delta}_M(\lambda,\Phi(\lambda,\cdot),t)$$

and (for $\lambda\to\infty$ with $\operatorname{Im}\lambda\le 0$)

$$\tilde{\hat{f}}(\lambda,\Phi(\lambda,\cdot),t)=\mathrm{O}(1/\lambda^3)\ \text{ and }\ \tilde{\delta}_M(\lambda,\Phi(\lambda,\cdot),t)=\mathrm{O}(1/\lambda^{M+3}). \tag{5.2}$$

Although (5.5) holds for $\lambda\in\mathbb{C}-\{0\}$ with $\operatorname{Im}\lambda\le 0$, the integral $\tilde{h}_M(\lambda,\omega(\lambda),t)$, is defined for every $\lambda\in\mathbb{C}-\{0\}$ with $\lambda\notin\{\pm i/\sqrt{\alpha(\tau)}:\,0\le\tau\le t\}$.

Now we proceed with the proof of Theorem 5 and we split it in several steps.

***Step 1*** Let us fix a $T>0$. In view of (1.5), for $x>0$ and $0<t<T$,

$$\int_{\lambda=-\infty}^{\infty}e^{i\lambda x-\Phi(\lambda,t)}\hat{u}_0(\lambda)d\lambda=\int_{-1}^{1}e^{i\lambda x-\Phi(\lambda,t)}\hat{u}_0(\lambda)d\lambda \tag{5.3}$$

$$+\left(\int_{-\infty}^{-1}+\int_{1}^{\infty}\right)e^{i\lambda x-\Phi(\lambda,t)}\sigma_M(\lambda)d\lambda+\left(\int_{-\infty}^{-1}+\int_{1}^{\infty}\right)e^{i\lambda x-\Phi(\lambda,t)}\upsilon_M(\lambda)d\lambda$$

$$=\int_{-1}^{1}e^{i\lambda x-\Phi(\lambda,t)}\hat{u}_0(\lambda)d\lambda+\sum_{k=1}^{M}u_0^{(k-1)}(0)\int_{\mathcal{K}_T}e^{i\lambda x-\Phi(\lambda,t)}\frac{d\lambda}{(i\lambda)^k}+\left(\int_{-\infty}^{-1}+\int_{1}^{\infty}\right)e^{i\lambda x-\Phi(\lambda,t)}(u_0^{(M)})\hat{\ }(\lambda)\frac{d\lambda}{(i\lambda)^M},$$

where $\mathcal{K}_T$ is a contour, starting at the point 1, going around the line segment $[c_0 i, c_T i]$ (in the upper half-plane) and ending at the point $-1$, as depicted in fig. 6, below.

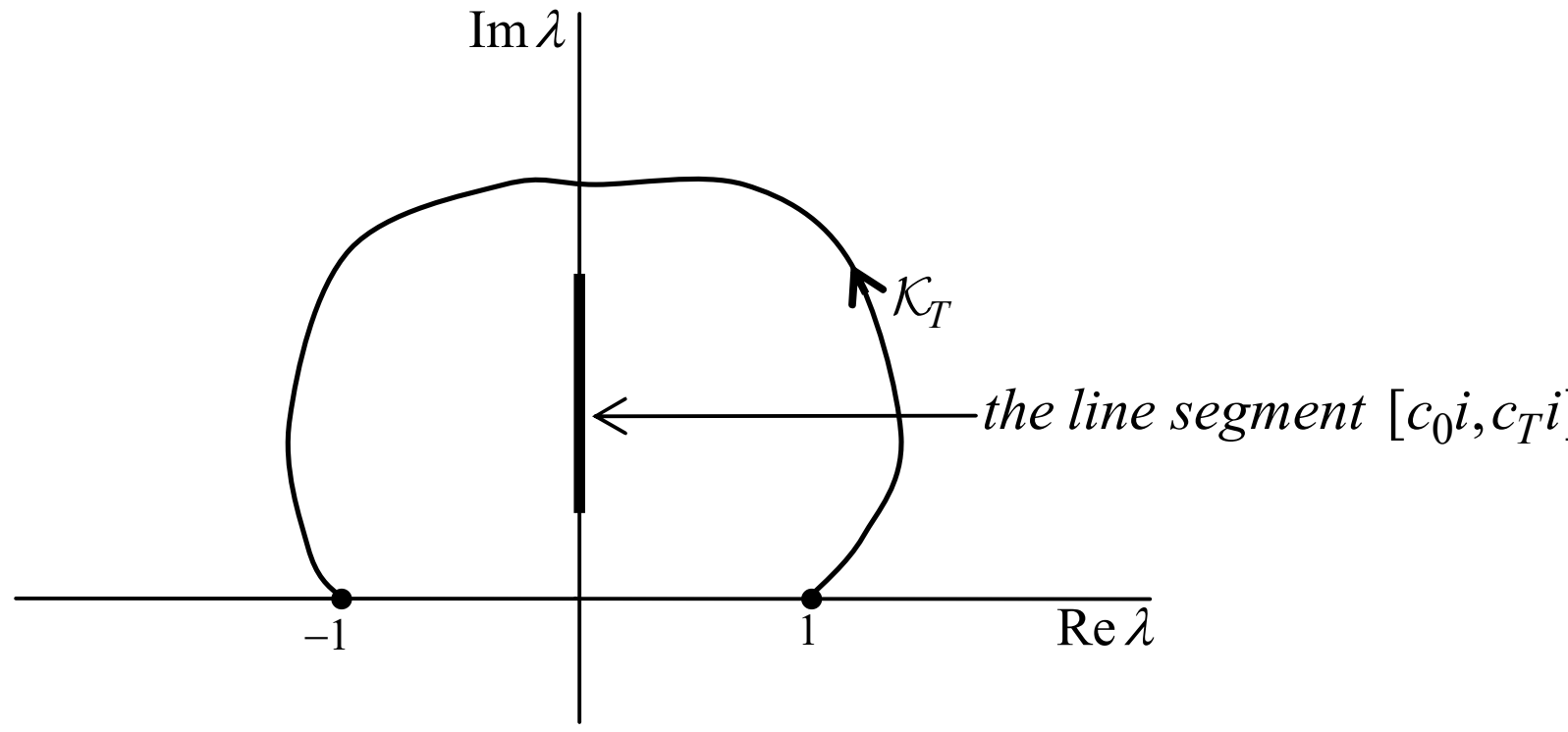

**Fig. 6** A typical contour $\mathcal{K}_T$

Since $M$ in (5.3) can be chosen arbitrarily large, we easily see that the function $\mathfrak{U}^+(x,t)$, defined by the integral in the LHS of (5.3), is $C^\infty$ for $x>0$ and $0<t<T$. It is at this stage that we justify that the integral $\mathfrak{U}^+(x,t)$ exists in the generalised sense, i.e., as limit of the form $\lim_{A\to\infty}\int_{-A}^{A}\cdots d\lambda$.

***Step 2*** Working as the previous step, we see that, in view of (5.1), the integral defining $\mathfrak{F}^+(x,t)$ can be written as follows: For $x>0$ and $0<t<T$,

$$(5.4)\quad \mathfrak{F}^+(x,t)=\int_{-1}^{1} e^{i\lambda x-\Phi(\lambda,t)}\widetilde{\hat{f}}(\lambda,\Phi(\lambda,\cdot),t)d\lambda+\int_{\mathcal{K}_T} e^{i\lambda x-\Phi(\lambda,t)}\widetilde{h}_M(\lambda,\Phi(\lambda,\cdot),t)d\lambda$$
$$+\left(\int_{-\infty}^{-1}+\int_{1}^{\infty}\right)e^{i\lambda x-\Phi(\lambda,t)}\widetilde{\delta}_M(\lambda,\Phi(\lambda,\cdot),t)d\lambda\,,$$

which, by the second relation of (5.2), implies that $\mathfrak{F}^+(x,t)$ is $C^\infty$ for $x>0$ and $0<t<T$. (We point out that, by the first relation of (5.2), that the integral defining $\mathfrak{F}^+(x,t)$ converges absolutely.)

***Step 3*** Fixing $T>0$ and restricting $t$ to the interval $0<t<T$, we may replace the contour $\mathcal{L}$, in the integrals defining the functions $\mathfrak{U}^-(x,t), \mathcal{I}_{g_0,g_0'}(x,t), \mathfrak{F}^-(x,t)$, by the contour $\mathcal{C}_T$. This implies that the functions $\mathfrak{U}^-(x,t), \mathcal{I}_{g_0,g_0'}(x,t), \mathfrak{F}^-(x,t)$ are $C^\infty$ for $x>0$ and $0<t<T$. Thus, combining this with the conclusions of steps 1 and 2, we conclude that, indeed, the function $u(x,t)$, defined by (4.13), is $C^\infty$ for $(x,t)\in Q$.

***Step 4*** Now we prove that the function $u(x,t)$ satisfies the differential equation of *Problem 2*, i.e., $\mathfrak{P}u=f$ in $Q$. Thinking of the action of the differential operator $\mathfrak{P}$ on the exponential part of the integrals of (4.13), namely $e^{i\lambda x-\Omega(\lambda,t)}$, and in particular in view of (5.4), and the computations made in *steps 1, 2* and *3*, it remains to show the following: For $x>0$ and $0<t<T$,

$$(5.5)\quad \int_{\mathcal{C}_T} i\lambda\left\{e^{i\lambda x-\Omega(\lambda,t)}\frac{\partial}{\partial t}[\mathfrak{G}(\lambda,t)]-\alpha(t)\frac{\partial}{\partial t}[\mathfrak{G}(\lambda,t)]\frac{\partial^2}{\partial x^2}[e^{i\lambda x-\Omega(\lambda,t)}]\right\}d\lambda=0\,,$$

$$(5.6)\quad \int_{-\infty}^{\infty}\left\{e^{i\lambda x-\Phi(\lambda,t)}\frac{\partial}{\partial t}\left[\widetilde{\hat{f}}(\lambda,\Phi(\lambda,\cdot),t)\right]-\alpha(t)\frac{\partial}{\partial t}\left[\widetilde{\hat{f}}(\lambda,\Phi(\lambda,\cdot),t)\right]\frac{\partial^2}{\partial x^2}[e^{i\lambda x-\Phi(\lambda,t)}]\right\}d\lambda=2\pi f(x,t)\,,$$

(5.7) $$\int_{\mathcal{C}_T}\left\{e^{i\lambda x-\Phi(\lambda,t)}\frac{\partial}{\partial t}\left[\tilde{\tilde{f}}(-\lambda,\Phi(\lambda,\cdot),t)\right]-\alpha(t)\frac{\partial}{\partial t}\left[\tilde{\tilde{f}}(-\lambda,\Phi(\lambda,\cdot),t)\right]\frac{\partial^2}{\partial x^2}[e^{i\lambda x-\Phi(\lambda,t)}]\right\}d\lambda=0.$$

Proof of (5.5) A straightforward computation shows that the integrand in (5.5) is equal to $i\lambda e^{i\lambda x}$. Thus, (5.5) follows from the equation $\int_{\mathcal{C}_T} i\lambda e^{i\lambda x}d\lambda=0$.

Proof of (5.6) The integrand in (5.6) is equal to $e^{i\lambda x}\hat{f}(\lambda,x)$, hence (5.6) follows from Fourier's inversion formula. (The integral in (5.6) is to be understood in the generalized sense and the previous computation relies, in part, also on (5.3) with $M=1$.)

Proof of (5.7) The integrand in (5.7) is equal to $e^{i\lambda x}\hat{f}(-\lambda,x)$, hence (5.7) follows from the equation $\int_{\mathcal{C}_T} e^{i\lambda x}\hat{f}(-\lambda,x)=0$.

***Step 5*** *Proof of (4.18)* Since $\mathfrak{G}(\lambda,t)\big|_{t=0}=0$, it follows from (4.15) that

(5.8) $$\lim_{t\to 0^+}\mathcal{I}_{g_0,g_0'}(x,t)=0.$$

Similarly, it is easy to see that

(5.9) $$\lim_{t\to 0^+}\mathfrak{F}^+(x,t)=0,\ \lim_{t\to 0^+}\mathfrak{F}^-(x,t)=0.$$

(To justify the first equation in (5.9), we use also the first relation of (5.2).)

Also

(5.10) $$\lim_{t\to 0^+}\int_{\mathcal{C}_T} e^{i\lambda x-\omega t}\hat{u}_0(-\lambda)d\lambda=\int_{\mathcal{C}_T} e^{i\lambda x}\hat{u}_0(-\lambda)d\lambda=0,$$

where the last equation follows from Cauchy's theorem.

On the other hand, it follows from (1.5) (with $M=1$) that

(5.11) $$\mathfrak{U}^+(x,t)=\int_{-1}^{1}e^{i\lambda x-\Phi(\lambda,t)}\hat{u}_0(\lambda)d\lambda+\int_{\mathcal{K}_T}e^{i\lambda x-\Phi(\lambda,t)}\frac{u_0(0)}{i\lambda}d\lambda+\left(\int_{-\infty}^{-1}+\int_{1}^{\infty}\right)e^{i\lambda x-\Phi(\lambda,t)}\frac{(u_0')\hat{}(\lambda)}{i\lambda}d\lambda.$$

Letting $t\to 0^+$, in (5.11), and using Cauchy's theorem and Jordan's lemma, we obtain

(5.12) $$\begin{aligned}\lim_{t\to 0^+}\mathfrak{U}^+(x,t)&=\int_{-1}^{1}e^{i\lambda x}\hat{u}_0(\lambda)d\lambda+\int_{\mathcal{K}_T}e^{i\lambda x}\frac{u_0(0)}{i\lambda}d\lambda+\left(\int_{-\infty}^{-1}+\int_{1}^{\infty}\right)e^{i\lambda x}\frac{(u_0')\hat{}(\lambda)}{i\lambda}d\lambda\\&=\int_{-1}^{1}e^{i\lambda x}\hat{u}_0(\lambda)d\lambda+\left(\int_{-\infty}^{-1}+\int_{1}^{\infty}\right)e^{i\lambda x}\frac{u_0(0)}{i\lambda}d\lambda+\left(\int_{-\infty}^{-1}+\int_{1}^{\infty}\right)e^{i\lambda x}\frac{(u_0')\hat{}(\lambda)}{i\lambda}d\lambda\\&=\int_{\lambda=-\infty}^{\infty}e^{i\lambda x}\hat{u}_0(\lambda)d\lambda=2\pi u_0(x).\end{aligned}$$

Thus, in view of (5.8) – (5.12), we obtain (4.18).

***Step 6*** *Proof of (4.19)* Firstly, we claim that

(5.13)
$$\lim_{x\to 0^+} \mathcal{I}_{g_0,g_0'}(x,t) = -\pi g_0(t) - \alpha(0) g_0(0) \int_{\mathcal{C}_T} \frac{i\lambda e^{-\Phi(\lambda,t)} d\lambda}{1+\alpha(0)\lambda^2}.$$

Indeed, by (4.17), we have

(5.14)
$$\lim_{x\to 0^+} \int_{\mathcal{C}_T} \frac{i\lambda e^{i\lambda x} d\lambda}{1+\alpha(t)\lambda^2} = -\frac{\pi}{\alpha(t)}.$$

On the other hand,

(5.15)
$$\lim_{x\to 0^+} \int_{\mathcal{C}_T} \left\{ i\lambda e^{i\lambda x - \Phi(\lambda,t)} \int_{\tau=0}^{t} e^{\Phi(\lambda,\tau)} \frac{\beta(\tau)-\alpha(\tau)\gamma(\tau)}{[1+\alpha(\tau)\lambda^2]^2} g_0(\tau)\, d\tau \right\} d\lambda$$
$$= \int_{\mathcal{C}_T} \left\{ i\lambda e^{-\Phi(\lambda,t)} \int_{\tau=0}^{t} e^{\Phi(\lambda,\tau)} \frac{\beta(\tau)+\alpha(\tau)\gamma(\tau)}{[1+\alpha(\tau)\lambda^2]^2} g_0(\tau)\, d\tau \right\} d\lambda$$
$$= \int_{-\infty}^{\infty} \left\{ i\lambda e^{-\Phi(\lambda,t)} \int_{\tau=0}^{t} e^{\Phi(\lambda,\tau)} \frac{\beta(\tau)+\alpha(\tau)\gamma(\tau)}{[1+\alpha(\tau)\lambda^2]^2} g_0(\tau)\, d\tau \right\} d\lambda = 0,$$

since the function in the $d\lambda$-integral is odd. (In the above computation, the equality $\int_{\mathcal{C}_T} \cdots d\lambda = \int_{-\infty}^{\infty} \cdots d\lambda$ follows from Cauchy's theorem and Jordan's lemma, and, for this to work, it is crucial that the integrand is $\mathrm{O}(1/\lambda^2)$. As a matter of fact, in this case, the integrand is $\mathrm{O}(1/\lambda^3)$.)
Similarly to (5.15), the limit, as $x \to 0^+$, of the last two terms in the RHS of (4.16), is also zero.
Thus, (5.13) follows from (4.16) and (5.14).
Next, we claim that

(5.16)
$$\lim_{x\to 0^+} [\mathfrak{U}^+(x,t) - \mathfrak{U}^-(x,t)] = 2u_0(0) \int_{\mathcal{K}_T} e^{-\Phi(\lambda,t)} \frac{d\lambda}{i\lambda}.$$

Indeed, in view of (5.1),
$$\mathfrak{U}^+(x,t) = \int_{-1}^{1} e^{i\lambda x - \Phi(\lambda,t)} \hat{u}_0(\lambda) d\lambda + \int_{\mathcal{K}_T} e^{i\lambda x - \Phi(\lambda,t)} \frac{u_0(0)}{i\lambda} d\lambda + \left( \int_{-\infty}^{-1} + \int_{1}^{\infty} \right) e^{i\lambda x - \Phi(\lambda,t)} \frac{(u_0')\hat{}(\lambda)}{i\lambda} d\lambda,$$

hence

(5.17)
$$\lim_{x\to 0^+} \mathfrak{U}^+(x,t) = \int_{-1}^{1} e^{-\Phi(\lambda,t)} \hat{u}_0(\lambda) d\lambda + \int_{\mathcal{K}_T} e^{-\Phi(\lambda,t)} \frac{u_0(0)}{i\lambda} d\lambda + \left( \int_{-\infty}^{-1} + \int_{1}^{\infty} \right) e^{-\Phi(\lambda,t)} \frac{(u_0')\hat{}(\lambda)}{i\lambda} d\lambda.$$

Similarly,

(5.18)
$$\lim_{x\to 0^+} \mathfrak{U}^-(x,t) = \int_{-1}^{1} e^{-\Phi(\lambda,t)} \hat{u}_0(-\lambda) d\lambda - \int_{\mathcal{K}_T} e^{-\Phi(\lambda,t)} \frac{u_0(0)}{i\lambda} d\lambda - \left( \int_{-\infty}^{-1} + \int_{1}^{\infty} \right) e^{-\Phi(\lambda,t)} \frac{(u_0')\hat{}(-\lambda)}{i\lambda} d\lambda.$$

Since
$$\int_{-1}^{1} e^{-\Phi(\lambda,t)} \hat{u}_0(\lambda) d\lambda = \int_{-1}^{1} e^{-\Phi(\lambda,t)} \hat{u}_0(-\lambda) d\lambda$$

and

$$\left(\int_{-\infty}^{-1}+\int_{1}^{\infty}\right)e^{-\Phi(\lambda,t)}\frac{(u_0')\hat{}(\lambda)}{i\lambda}d\lambda=-\left(\int_{-\infty}^{-1}+\int_{1}^{\infty}\right)e^{-\Phi(\lambda,t)}\frac{(u_0')\hat{}(-\lambda)}{i\lambda}d\lambda,$$

equation (5.16) follows from (5.17) and (5.18).

Also, in view of (5.2),

$$\lim_{x\to 0^+}\mathfrak{F}^+(x,t)=\int_{-\infty}^{\infty}\left\{e^{-\Phi(\lambda,t)}\int_{\tau=0}^{t}\frac{e^{\Phi(\lambda,\tau)}}{1+\alpha(\tau)\lambda^2}\hat{f}(\lambda,\tau)\,d\tau\right\}d\lambda$$

and

$$\lim_{x\to 0^+}\mathfrak{F}^-(x,t)=\int_{-\infty}^{\infty}\left\{e^{-\Phi(\lambda,t)}\int_{\tau=0}^{t}\frac{e^{\Phi(\lambda,\tau)}}{1+\alpha(\tau)\lambda^2}\hat{f}(-\lambda,\tau)\,d\tau\right\}d\lambda,$$

and, therefore,

$$\lim_{x\to 0^+}[\mathfrak{F}^+(x,t)-\mathfrak{F}^-(x,t)]=0. \tag{5.19}$$

Finally, we claim that

$$\int_{\mathcal{K}_T}e^{-\Phi(\lambda,t)}\frac{d\lambda}{\lambda}=\alpha(0)\int_{\mathcal{C}_T}\frac{\lambda e^{-\Phi(\lambda,t)}d\lambda}{1+\alpha(0)\lambda^2}. \tag{5.20}$$

Since

$$\int_{-1}^{1}\frac{\lambda e^{-\Phi(\lambda,t)}d\lambda}{1+\alpha(0)\lambda^2}=0,$$

we have

$$\int_{\mathcal{C}_T}\frac{\lambda e^{-\Phi(\lambda,t)}d\lambda}{1+\alpha(0)\lambda^2}=\int_{\mathcal{K}_T}\frac{\lambda e^{-\Phi(\lambda,t)}d\lambda}{1+\alpha(0)\lambda^2}.$$

Also,

$$\int_{\mathcal{K}_T}e^{-\Phi(\lambda,t)}\frac{d\lambda}{\lambda}-\alpha(0)\int_{\mathcal{K}_T}\frac{\lambda e^{-\Phi(\lambda,t)}d\lambda}{1+\alpha(0)\lambda^2}=\int_{\mathcal{K}_T}\frac{e^{-\Phi(\lambda,t)}}{\lambda[1+\alpha(0)\lambda^2]}d\lambda=\left(\int_{-\infty}^{-1}+\int_{1}^{\infty}\right)\frac{e^{-\Phi(\lambda,t)}}{\lambda[1+\alpha(0)\lambda^2]}d\lambda=0,$$

and (5.20) follows.

Now, (4.19) follows from (4.14), (5.13), (5.16), (5.19), (5.20) and the assumption “$u_0(0)=g_0(0)$”.

***Completion of the proof*** The conclusions of the theorem follow from the computations and the results of *Steps* 1 – 6.

**Proof of Theorem 6** For the integral $\mathfrak{U}^+(x,t)$, differentiating (5.3), we obtain, for $x>0$ and $0<t<T$, that

$$\frac{\partial^{n+m}\mathfrak{U}^+(x,t)}{\partial x^n\partial t^m}=\int_{-1}^{1}(i\lambda)^n e^{i\lambda x}\frac{\partial^m}{\partial t^m}[e^{-\Phi(\lambda,t)}]\hat{u}_0(\lambda)d\lambda \tag{5.21}$$

$$+\sum_{k=1}^{M}u_0^{(k-1)}(0)\int_{\mathcal{K}_T}(i\lambda)^n e^{i\lambda x}\frac{\partial^m}{\partial t^m}[e^{-\Phi(\lambda,t)}]\frac{d\lambda}{(i\lambda)^k}+\left(\int_{-\infty}^{-1}+\int_{1}^{\infty}\right)(i\lambda)^n e^{i\lambda x}\frac{\partial^m}{\partial t^m}[e^{-\Phi(\lambda,t)}](u_0^{(M)})\hat{}(\lambda)\frac{d\lambda}{(i\lambda)^M}.$$

Similarly, for the integral $\mathfrak{F}^{+}(x,t)$, differentiating (5.4), we obtain

$$(5.22)\quad \frac{\partial^{n+m}\mathfrak{F}^{+}(x,t)}{\partial x^{n}\partial t^{m}} = \int_{-1}^{1}(i\lambda)^{n}e^{i\lambda x}\frac{\partial^{m}}{\partial t^{m}}[e^{-\Phi(\lambda,t)}\tilde{\tilde{f}}(\lambda,\Phi(\lambda,\cdot),t)]d\lambda$$

$$+\int_{\mathcal{K}_T}(i\lambda)^{n}e^{i\lambda x}\frac{\partial^{m}}{\partial t^{m}}[e^{-\Phi(\lambda,t)}\tilde{h}_M(\lambda,\Phi(\lambda,t),\cdot)]d\lambda$$

$$+\left(\int_{-\infty}^{-1}+\int_{1}^{\infty}\right)(i\lambda)^{n}e^{i\lambda x}\frac{\partial^{m}}{\partial t^{m}}[e^{-\Phi(\lambda,t)}\tilde{\delta}_M(\lambda,\Phi(\lambda,\cdot),t)]d\lambda .$$

Recalling that the contour $\mathcal{L}$ in the integrals $\mathfrak{U}^{-}(x,t)$, $\mathcal{I}_{g_0,g_0'}(x,t)$ and $\mathfrak{F}^{-}(x,t)$, of (4.13), can be replaced by the finite contour $\mathcal{C}_T$ (when $t$ is restricted to the interval $0<t<T$), we see that we can immediately interchange the differentiations $\partial^{n+m}/\partial x^{n}\partial t^{m}$ with the integrations, obtaining formulas analogous to (5.21) and (5.22).

Now the 1st and 2nd conclusions of the theorem follow from these formulas.

Having proved the 1st and 2nd assertions, the 3rd conclusion follows. Combining this with the conclusions of Theorem 1, we easily obtain the 4th and 5th conclusions.

Finally, for the proof of the 6th assertion, it suffices to appeal once more to the formulas (5.21), (5.22) and their analogues, mentioned above. This completes the proof of the theorem.

**Proof of Theorem 7** It is similar to the proof of the analogous theorem for the autonomous case discussed in [17].

## 5. Applications

The solutions and their qualitative properties obtained above are directly useful for further analysis:

**Further results for the Barenblatt equation**

**5.1** ***The global relation*** By Theorems 5, 6 and 7 (assuming (1.3), (1.4), (1.9), (1.10) and "$u_0(0)=g_0(0)$"), the derivation of (4.8) is rigorous, i.e., the function $u(x,t)$, defined by (4.13), satisfies, indeed, these equations. The point here is that for the derivation of the GR (4.8), we had assumed that there existed a solution of Problem 2 such that the various steps in the process of the derivation were justified. Now, having proved the Theorems 5, 6, and 7, we know that such a solution exists, making the derivation – a posteriori – rigorous and the GR (4.8) valid. This is also essential for:

**5.2** ***A non-controlability theorem*** *Suppose that $u_0(x)\in\mathcal{S}([0,\infty))$, $\hat{u}_0(\lambda)$ does not extend to an analytic function of $\lambda\in\mathbb{C}-\{[ic_0,\infty i)\cup[-ic_0,-\infty i)\}$ and $f$ satisfies (1.4). Then for every $t_0>0$, so that $\tilde{\tilde{f}}(\lambda,\Phi(\lambda,\cdot),t_0)$ extends to an analytic of $\lambda\in\mathbb{C}-\{[ic_0,\infty i)\cup[-ic_0,-\infty i)\}$, and for every $g_0$ with $g_0(0)=u_0(0)$, we have that the UTM solution $u(x,t)$ of Problem 2, given by (4.13), satisfies $u(x,t_0)\not\equiv 0$.* (See [49] for connections to *natural* boundaries.)

**Proof** As we pointed out above, the GR given by (4.8) holds valid. Now, if $u(x,t_0)\equiv 0$, for some $t_0>0$, then $\hat{u}(\lambda,t_0)\equiv 0$. But this would contradict (4.8) when $t=t_0$, since all the terms in $(4.8)\big|_{t=t_0}$, except $\hat{u}_0(\lambda)$, are analytic function of $\lambda\in\mathbb{C}-\{[ic_0,\infty i)\cup[-ic_0,-\infty i)\}$.

*Remark* We point out that it is because of Theorems 5 and 6, that the integral in (4.8), which involves the functions

(5.1) $$g_1'(\tau)=\frac{d}{d\tau}[u_x(0,\tau)],\ g_1(\tau)=u_x(0,\tau),$$

define analytic functions for $\lambda\in\mathbb{C}-\{[ic_0,\infty i)\cup[-ic_0,-\infty i)\}$. Indeed, in order to ensure this analyticity, we need that the functions (5.1) be smooth for $\tau\in[0,\infty)$.

**5.3** ***A uniqueness theorem*** *Let us assume, in addition to the assumptions of Theorem 6, that* $\alpha'(t)\le 2\beta(t)$ *and* $\gamma(t)\le 0$ $(\forall t>0)$. *Then the solution* $u(x,t)$ *of Problem 2, given by (4.13), is unique in the following sense: If* $v(x,t)$ *is* $C^3$ *in* $\overline{Q}-\{(0,0)\}$, *satisfies (1.2),*

*(5.2)* $$\lim_{x\to\infty} v(x,t)=0,\ \sup_{x\ge1}|v_x(x,t)|<\infty,\ \sup_{x\ge1}|v_{xt}(x,t)|<\infty\ (\forall t>0),$$

*and, for every* $T>0$, *the functions* $|v(x,t)|^2$, $|v_x(x,t)|^2$, $|v_t(x,t)|^2$, $|v_{xx}(x,t)|^2$, $|v_{xt}(x,t)|^2$, *are, uniformly for* $0<t\le T$, *integrable with respect to* $x\in[0,\infty)$, *i.e., there exists a positive function* $\mathcal{M}_T(x)$ *such that* $\int_0^\infty \mathcal{M}_T(x)dx<+\infty$ *and, for* $0<t\le T$ *and* $x>0$,

*(5.3)* $$\max\{|v(x,t)|^2,|v_x(x,t)|^2,|v_t(x,t)|^2,|v_{xx}(x,t)|^2,|v_{xt}(x,t)|^2\}\le\mathcal{M}_T(x),$$

*then* $v\equiv u$.

**Proof** By Theorems 5, 6 and 7, and the assumptions on $v$, the difference

$$w(x,t):=v(x,t)-u(x,t)$$

satisfies the following:

(5.4) $$w(x,0)=0\ \text{for}\ x>0,\ w_x(x,0)=0\ \text{for}\ x>0,\ w(0,t)=0\ \text{for}\ t>0,$$

and

(5.5) $$w_t(x,t)-\alpha(t)w_{xxt}(x,t)-\beta(t)w_{xx}(x,t)-\gamma(t)w(x,t)=0,$$

and (5.2), (5.3) hold with $v$ replaced by $w$ (by taking $\mathcal{M}_T$ sufficiently large – if necessary). This last conclusion follows from the fact that the function $u$ is rapidly decreasing as $x\to\infty$, uniformly for $0<t\le T$, by Theorem 7, and that $u$ and its derivatives are bounded near the point $(0,0)$, by Theorem 6. Therefore, the functions $ww_t$, $ww_{xxt}$ and $ww_{xx}$ are absolutely integrable, with respect to $x$, in $[0,+\infty)$, and, (5.5) implies that, for fixed $t>0$,

(5.6) $$\int_0^\infty [w(x,t)w_t(x,t)-\alpha(t)w(x,t)w_{xxt}(x,t)-\beta(t)w(x,t)w_{xx}(x,t)-\gamma(t)w(x,t)w(x,t)]dx=0.$$

By Lebesgue's dominated convergence theorem, (5.6) can be written as follows:

(5.7) $$\frac{d}{dt}\left(\int_0^\infty \tfrac{1}{2}[w(x,t)]^2dx\right)-\alpha(t)[w(x,t)w_{xt}(x,t)]\Big|_{x=0}^{\infty}+\alpha(t)\int_0^\infty w_x(x,t)w_{xt}(x,t)dx$$
$$-\beta(t)[w(x,t)w_x(x,t)]\Big|_{x=0}^{\infty}+\beta(t)\int_0^\infty [w_x(x,t)]^2dx-\gamma(t)\int_0^\infty [w(x,t)]^2dx=0.$$

But, by (5.2) (applied with $w$ in place of $v$) and (5.4),

$$[w(x,t)w_{xt}(x,t)]\big|_{x=0}^{\infty} = \lim_{A\to\infty} w(A,t)w_{xt}(A,t) = 0 \text{ and } [w(x,t)w_x(x,t)]\big|_{x=0}^{\infty} = \lim_{A\to\infty} w(A,t)w_x(A,t) = 0.$$

Thus, (5.7) becomes

$$\frac{d}{dt}\left(\int_0^\infty \tfrac{1}{2}[w(x,t)]^2 dx\right) + \alpha(t)\int_0^\infty w_x(x,t)w_{xt}(x,t)dx + \beta(t)\int_0^\infty [w_x(x,t)]^2 dx - \gamma(t)\int_0^\infty [w(x,t)]^2 dx = 0.$$

Therefore,

$$\frac{d}{dt}\left(\int_0^\infty \{\tfrac{1}{2}[w(x,t)]^2 + \tfrac{\alpha(t)}{2}[w_x(x,t)]^2\}dx\right)$$
$$= \int_0^\infty \tfrac{\alpha'(t)}{2}[w_x(x,t)]^2\}dx - \beta(t)\int_0^\infty [w_x(x,t)]^2 dx + \gamma(t)\int_0^\infty [w(x,t)]^2 dx \le 0,$$

where used also the assumptions "$\alpha'(t) \le 2\beta(t)$" and "$\gamma(t) \le 0$".

It follows that the function $\int_0^\infty \{[w(x,t)]^2 + \alpha(t)[w_x(x,t)]^2\}dx$ is non-increasing for $t \in (0,+\infty)$, i.e.,

$$\int_0^\infty \{[w(x,t)]^2 + \alpha(t)[w_x(x,t)]^2\}dx \le \int_0^\infty \{[w(x,s)]^2 + \alpha(s)[w_x(x,s)]^2\}dx \text{ for } 0 < s < t.$$

Letting $s \to 0^+$, we obtain

$$\int_0^\infty \{[w(x,t)]^2 + \alpha(t)[w_x(x,t)]^2\}dx \le \int_0^\infty \{[w(x,0)]^2 + \alpha(0)[w_x(x,0)]^2\}dx = 0, \text{ for } t > 0.$$

Therefore, $[w(x,t)]^2 + \alpha(t)[w_x(x,t)]^2 = 0$, for $t > 0$ and $x > 0$, hence $w(x,t) \equiv 0$.

**Further results for the heat equation**

Similarly to the above, we can prove analogous results for the solution of Problem 1.

**5.4** ***A non-controlability theorem*** *Suppose that $u_0(x) \in \mathcal{S}([0,\infty))$, $\hat{u}_0(\lambda)$ does not extend to an analytic function of $\lambda \in \mathbb{C}$ and $f$ satisfies (1.4). Then for every $t_0 > 0$, so that $\tilde{\hat{f}}(\lambda,\Omega(\lambda,\cdot),t_0)$ extends to an analytic of $\lambda \in \mathbb{C}$, and for every $g_0$ for which*

$$g_0(0) = u_0(0) \text{ and } \left.\frac{dg_0(t)}{dt}\right|_{t=0} = \alpha(0)\left.\frac{d^2u_0(x)}{dx^2}\right|_{x=0} + \beta(0)g_0(0) + f(0,0),$$

*we have that the UTM solution $u(x,t)$ of Problem 1, given by (2.9), satisfies $u(x,t_0) \not\equiv 0$.*

Remark For rigorous UTM-based controllability results in autonomous cases see e.g. [28, 12].

**5.5** ***A uniqueness theorem*** *The function $u(x,t)$, defined by (2.9), is the unique solution of Problem 1, in the following sense: If we assume (1.3), (1.4), (2.17) and $\beta(t) \le 0$ ($\forall t > 0$),*

*$v(x,t)$ is a $C^2$ function in $\overline{Q} - \{(0,0\}$ and satisfies (1.1),*

*$\lim_{x\to\infty} v(x,t) = 0$ and $\sup_{x\ge1}|v_x(x,t)| < \infty$ ($\forall t > 0$),*

*and, for every $T > 0$,*

$|v(x,t)|^2$, $|v_t(x,t)|^2$ *are, uniformly for* $0<t\le T$, *integrable with respect to* $x\in[0,\infty)$,
*i.e., there exists a positive function* $\mathcal{M}_T(x)$ *such that* $\int_0^\infty \mathcal{M}_T(x)dx<+\infty$ *and, for* $0<t\le T$,

$$|v(x,t)|^2 \le \mathcal{M}_T(x) \text{ and } |v_t(x,t)|^2 \le \mathcal{M}_T(x) \ (x>0),$$

*then* $v\equiv u$.

<u>*Remark*</u> For rigorous UTM-based uniqueness results in autonomous cases see e.g. [28,16-18,22,23].

**An eventual periodicity theorem for the Barenblantt equation**

**5.6** ***Theorem*** *Let us consider Problem 2, with the function* $\alpha(t)\equiv\alpha$, *being constant, and* $\gamma(t)\equiv 0$, *in which case the problem becomes:*

$$\text{(5.8)}\qquad \begin{cases}\partial_t u=\alpha\partial_{xxt}u+\beta(t)\partial_{xx}u+f(x,t), & (x,t)\in Q=\mathbb{R}^+\times\mathbb{R}^+,\\ u(x,0)=u_0(x), & x\in\mathbb{R}^+,\\ u(0,t)=g_0(t), & t\in\mathbb{R}^+.\end{cases}$$

*Furthermore we assume that for a constant* $T>0$,

$$\text{(5.9)}\qquad g_0(t+T)=g_0(t),\ f(x,t+T)=f(x,t) \text{ and } \beta(t+T)=\beta(t),\ \forall t,$$

*and that* $\beta(t)\ge 0$ *with* $B(t):=\int_0^t\beta(\tau)d\tau\to\infty$, *as* $t\to\infty$.

*Then, for the solution* $u(x,t)$ *of (5.8), we have: For fixed* $x>0$,

$$\lim_{t\to\infty}[u(x,t+T)-u(x,t)]=0.$$

*More precisely, for fixed* $x>0$ *and as* $t\to\infty$,

$$\text{(5.10)}\quad u(x,t+T)-u(x,t)=\frac{x}{2\sqrt{\pi}}\left(\int_{-T}^{0}\beta(\tau)g_0(\tau)d\tau+\int_{-T}^{0}\int_0^\infty yf(y,\tau)dyd\tau\right)\frac{1}{[B(t)]^{3/2}}+\mathrm{O}(1/[B(t)]^{5/2}).$$

*Examples of functions* $\beta(t)$: 1. $\beta(t)=3+\cos t-2\sin(2t)$, $B(t)=3t+\sin t+\cos(2t)$

2. $\beta(t)=e^{\sin t}$, $B(t)=\int_0^t e^{\sin\tau}d\tau$.

*Remark* More generally, if $\beta(t)$ is any non-negative and $T-$periodic function, and its mean value in the interval $[0,T]$

$$\mu:=\frac{1}{T}\int_0^T\beta(\tau)d\tau$$

is a positive number, then

$$B(t)=\int_0^t\beta(\tau)d\tau=\mu t+\mathrm{O}(1),$$

and, therefore, $B(t)\to\infty$, as $t\to\infty$. Indeed, if $t=N(t)T+\varepsilon(t)$, where $N(t)\in\mathbb{N}$ and $0\le\varepsilon(t)<T$, then

$$\frac{1}{T}\int_0^t\beta(\tau)d\tau-\mu t=[N(t)-t]\mu+\frac{1}{T}\int_0^{\varepsilon(t)}\beta(\tau)d\tau,$$

and this implies that

$$\left|\frac{1}{T}\int_0^t \beta(\tau)d\tau - \mu t\right| \le \mu .$$

**Proof of the theorem** Setting

$$\omega(\lambda) = -\frac{\lambda^2}{1+\alpha\lambda^2},$$

we have

$$\phi(\lambda,t) = \beta(t)\omega(\lambda) \text{ and } \Phi(\lambda,t) = \int_{\tau=0}^{t} \phi(\lambda,\tau)d\tau = -B(t)\omega(\lambda),$$

and the solution $u(x,t)$ of (5.8) takes the form:

$$\begin{aligned}(5.11)\quad 2\pi u(x,t) &= \int_{-\infty}^{\infty} e^{i\lambda x+B(t)\omega(\lambda)}\hat{u}_0(\lambda)d\lambda - \int_{\mathcal{C}} e^{i\lambda x+B(t)\omega(\lambda)}\hat{u}_0(-\lambda)d\lambda \\ &\quad -2\int_{\mathcal{C}}\left\{\frac{i\lambda}{1+\alpha\lambda^2}e^{i\lambda x+B(t)\omega(\lambda)}\int_{\tau=0}^{t} e^{-B(\tau)\omega(\lambda)}[\alpha g_0'(\tau)+\beta(\tau)g_0(\tau)]\,d\tau\right\}d\lambda \\ &\quad + \int_{-\infty}^{\infty}\left\{e^{i\lambda x+B(t)\omega(\lambda)}\int_{\tau=0}^{t}\frac{e^{-B(\tau)\omega(\lambda)}}{1+\alpha\lambda^2}\hat{f}(\lambda,\tau)\,d\tau\right\}d\lambda - \int_{\mathcal{C}}\left\{e^{i\lambda x+B(t)\omega(\lambda)}\int_{\tau=0}^{t}\frac{e^{-B(\tau)\omega(\lambda)}}{1+\alpha\lambda^2}\hat{f}(-\lambda,\tau)\,d\tau\right\}d\lambda \\ &=: I_1(x,t) - I_2(x,t) - 2I_3(x,t) + I_4(x,t) - I_5(x,t).\end{aligned}$$

Throughtout this proof we should keep in mind the following remark:

If $\lambda = re^{i\theta}$, with $r = |\lambda|$ and $\theta = \arg\lambda$, we have

$$\left|e^{B(t)\omega(\lambda)}\right| = \exp\left(-B(t)r^2\frac{\alpha r^2+\cos 2\theta}{\left|1+\alpha r^2 e^{i2\theta}\right|^2}\right).$$

In particular, if $\lambda\in\mathbb{C}$ and $|\lambda|\ge 1/\sqrt{\alpha}$, then $\left|e^{B(t)\omega(\lambda)}\right|\le 1$. Also, $\left|e^{B(t)\omega(\lambda)}\right|\le 1$, for every $\lambda\in\mathbb{R}$.

*Step 1* Let us write the following version of (5.3):

$$\begin{aligned}(5.12)\quad I_1(x,t) &= \int_{-\infty}^{\infty} e^{i\lambda x+B(t)\omega(\lambda)}\hat{u}_0(\lambda)d\lambda = \int_{-1/\sqrt{\alpha}}^{1/\sqrt{\alpha}} e^{i\lambda x+B(t)\omega(\lambda)}\hat{u}_0(\lambda)d\lambda \\ &\quad + \left(\int_{-\infty}^{-1/\sqrt{\alpha}} + \int_{1/\sqrt{\alpha}}^{\infty}\right)e^{i\lambda x+B(t)\omega(\lambda)}\sigma_M(\lambda)d\lambda + \left(\int_{-\infty}^{-1/\sqrt{\alpha}} + \int_{1/\sqrt{\alpha}}^{\infty}\right)e^{i\lambda x+B(t)\omega(\lambda)}\upsilon_M(\lambda)d\lambda \\ &= \int_{-1/\sqrt{\alpha}}^{1/\sqrt{\alpha}} e^{i\lambda x+B(t)\omega(\lambda)}\hat{u}_0(\lambda)d\lambda + \sum_{k=1}^{M} u_0^{(k-1)}(0)\int_{\mathcal{K}_\alpha} e^{i\lambda x+B(t)\omega(\lambda)}\frac{d\lambda}{(i\lambda)^k} \\ &\quad + \left(\int_{-\infty}^{-1/\sqrt{\alpha}} + \int_{1/\sqrt{\alpha}}^{\infty}\right)e^{i\lambda x+B(t)\omega(\lambda)}(u_0^{(M)})\hat{}(\lambda)\frac{d\lambda}{(i\lambda)^M} \\ &=: I_{1,1}(x,t) + I_{1,2,M}(x,t) + I_{1,3,M}(x,t),\end{aligned}$$

where $\mathcal{K}_\alpha$ is a simple contour, which starts at the point $\lambda = 1/\sqrt{\alpha}$, ends at the point $\lambda = -1/\sqrt{\alpha}$ and lies in the set $\{\operatorname{Im}\lambda > 0, |\lambda| > 1/\sqrt{\alpha}\} \cup \{1/\sqrt{\alpha}, -1/\sqrt{\alpha}\}$. (See fig. 7.)

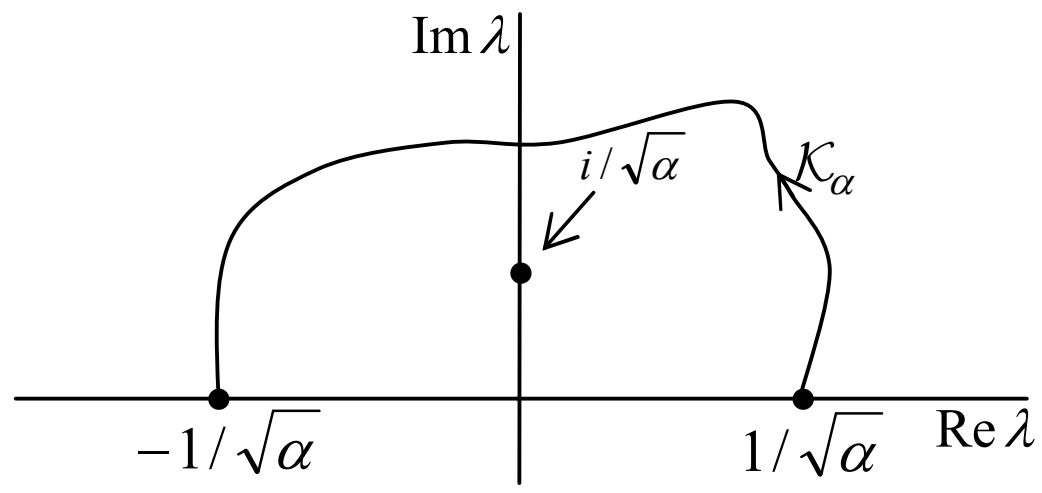

**Fig. 7** A typical contour $\mathcal{K}_\alpha$

Also,

(5.13) $$I_2(x,t) = \int_{\mathcal{C}} e^{i\lambda x + B(t)\omega(\lambda)} \hat{u}_0(-\lambda)d\lambda = \int_{-1/\sqrt{\alpha}}^{1/\sqrt{\alpha}} e^{B(t)\omega(\lambda)} e^{i\lambda x} \hat{u}_0(-\lambda)d\lambda + \int_{\mathcal{K}_\alpha} e^{B(t)\omega(\lambda)} e^{i\lambda x} \hat{u}_0(-\lambda)d\lambda.$$
$$=: I_{2,1}(x,t) + I_{2,2}(x,t).$$

Using Laplace asymptotic expansion formula, we obtain

$$I_{1,1}(x,t) - I_{2,1}(x,t) = \int_{-1/\sqrt{\alpha}}^{1/\sqrt{\alpha}} e^{B(t)\omega(\lambda)} e^{i\lambda x} [\hat{u}_0(\lambda) - \hat{u}_0(-\lambda)]d\lambda$$
$$= x\sqrt{\pi}\left(\int_0^\infty y u_0(y)dy\right)\frac{1}{[B(t)]^{3/2}} + \mathrm{O}(1/[B(t)]^{5/2}).$$

On the other hand, using the identity

$$e^{B(t)\omega(\lambda)} = \frac{1}{B(t)} \frac{1}{\omega'(\lambda)} \frac{d[e^{B(t)\omega(\lambda)}]}{d\lambda},$$

and integrating by parts, repeatedly, we see that

$$I_{1,2}(x,t) = \sum_{k=1}^{M} u_0^{(k-1)}(0) \int_{\mathcal{K}_\alpha} e^{i\lambda x + B(t)\omega(\lambda)} \frac{d\lambda}{(i\lambda)^k} = \mathrm{O}(1/[B(t)]^{\ell}), \text{ for every } \ell \in \mathbb{N}.$$

Next, choosing $M$ sufficiently large, each time, we show that the sum of the integrals

$$I_{1,2,M}(x,t) + I_{1,3,M}(x,t)$$

is also $\mathrm{O}(1/[B(t)]^{\ell})$, for every $\ell \in \mathbb{N}$. (Of course, the above sum is independent of $M$.)
Similarly, we see that the integral $I_{2,2}(x,t)$ is also $\mathrm{O}(1/[B(t)]^{\ell})$, for every $\ell \in \mathbb{N}$.
Therefore, in view of (5.12) and (5.13),

$$I_1(x,t) - I_2(x,t) = x\sqrt{\pi}\left(\int_0^\infty y u_0(y)dy\right)\frac{1}{[B(t)]^{3/2}} + \mathrm{O}(1/[B(t)]^{5/2}).$$

It follows that, for a fixed $x > 0$, as $t \to \infty$,

(5.14) $$[I_1(x,t+T) - I_2(x,t+T)] - [I_1(x,t) - I_2(x,t)] = \mathrm{O}(1/[B(t)]^{5/2}).$$

*Reminder – Laplace asymptotic expansion formula* The version of Laplace formula that we are using is the following: If $J\subset\mathbb{R}$ is a bounded open interval, $\lambda_0\in J$, $\omega'(\lambda_0)=0$, $\omega'(\lambda)\neq 0$ for $\lambda\in\bar{J}-\{\lambda_0\}$, and $\omega''(\lambda_0)<0$ then, as $t\to\infty$,

$$\int_J e^{B(t)\omega(\lambda)}h(\lambda)d\lambda=e^{t\omega(\lambda_0)}\frac{\sqrt{2\pi}}{\sqrt{-\omega''(\lambda_0)}}\left\{h(\lambda_0)\frac{1}{\sqrt{B(t)}}\right.$$
$$\left.+\left[-\frac{h^{(2)}(\lambda_0)}{2\omega^{(2)}(\lambda_0)}+\frac{h(\lambda_0)\omega^{(4)}(\lambda_0)}{8[\omega^{(2)}(\lambda_0)]^2}+\frac{h^{(1)}(\lambda_0)\omega^{(3)}(\lambda_0)}{2[\omega^{(2)}(\lambda_0)]^2}-\frac{5h(\lambda_0)[\omega^{(3)}(\lambda_0)]^2}{24[\omega^{(2)}(\lambda_0)]^3}\right]\frac{1}{[B(t)]^{3/2}}+\mathrm{O}(1/[B(t)]^{5/2})\right\}.$$

The functions $\omega$ and $h$ are assumed to be sufficiently smooth in $\bar{J}=closure(J)$.

*Note* For the application of Laplace formula, we note the following:

$$\omega(\lambda)=-\frac{\lambda^2}{1+\alpha\lambda^2},\ \omega'(\lambda)=-\frac{2\lambda}{(1+\alpha\lambda^2)^2},\ \omega''(\lambda)=-\frac{2}{(1+\alpha\lambda^2)^2}+\frac{8\alpha\lambda^2}{(1+\alpha\lambda^2)^3},$$
$$\omega'(\lambda)=0\ \Leftrightarrow\ \lambda=0,$$
$$\omega(0)=0,$$
$$\omega''(0)=-2,$$
$$\omega'''(0)=0.$$

*Step 2* We have

$$I_3(x,t+T)-I_3(x,t)=\int_C\left\{\frac{i\lambda}{1+\alpha\lambda^2}e^{i\lambda x+B(t+T)\omega(\lambda)}\int_{\tau=0}^{t+T}e^{-B(\tau)\omega(\lambda)}[\alpha g_0'(\tau)+\beta(\tau)g_0(\tau)]\,d\tau\right\}d\lambda$$
$$-\int_C\left\{\frac{i\lambda}{1+\alpha\lambda^2}e^{i\lambda x+B(t)\omega(\lambda)}\int_{\tau=0}^{t}e^{-B(\tau)\omega(\lambda)}[\alpha g_0'(\tau)+\beta(\tau)g_0(\tau)]\,d\tau\right\}d\lambda.$$

Now, in view of the periodicity assumptions (5.9),

$$e^{B(t+T)\omega(\lambda)}\int_{\tau=0}^{t+T}e^{-B(\tau)\omega(\lambda)}[\alpha g_0'(\tau)+\beta(\tau)g_0(\tau)]d\tau-e^{B(t)\omega(\lambda)}\int_{\tau=0}^{t}e^{-B(\tau)\omega(\lambda)}[\alpha g_0'(\tau)+\beta(\tau)g_0(\tau)]d\tau$$
$$=\exp\left(\omega(\lambda)\int_{\tau=0}^{t+T}\beta(\tau)d\tau\right)\int_{\tau=0}^{t+T}\exp\left(-\omega(\lambda)\int_{\sigma=0}^{\tau}\beta(\sigma)d\sigma\right)\beta(\tau)g_0(\tau)d\tau$$
$$-\exp\left(\omega(\lambda)\int_{\tau=0}^{t}\beta(\tau)d\tau\right)\int_{\tau=0}^{t}\exp\left(-\omega(\lambda)\int_{\sigma=0}^{\tau}\beta(\sigma)d\sigma\right)[\alpha g_0'(\tau)+\beta(\tau)g_0(\tau)]d\tau$$
$$=\exp\left(\omega(\lambda)\int_{\tau=-T}^{t}\beta(\tau)d\tau\right)\int_{\tau=-T}^{t}\exp\left(-\omega(\lambda)\int_{\sigma=-T}^{\tau}\beta(\sigma)d\sigma\right)[\alpha g_0'(\tau)+\beta(\tau)g_0(\tau)]d\tau$$
$$-\exp\left(\omega(\lambda)\int_{\tau=0}^{t}\beta(\tau)d\tau\right)\int_{\tau=0}^{t}\exp\left(-\omega(\lambda)\int_{\sigma=0}^{\tau}\beta(\sigma)d\sigma\right)[\alpha g_0'(\tau)+\beta(\tau)g_0(\tau)]d\tau$$
$$=\exp\left(\omega(\lambda)\int_{\tau=0}^{t}\beta(\tau)d\tau\right)\int_{\tau=-T}^{t}\exp\left(-\omega(\lambda)\int_{\sigma=0}^{\tau}\beta(\sigma)d\sigma\right)[\alpha g_0'(\tau)+\beta(\tau)g_0(\tau)]d\tau$$
$$\exp\left(\omega(\lambda)\int_{\tau=0}^{t}\beta(\tau)d\tau\right)\int_{\tau=0}^{t}\exp\left(-\omega(\lambda)\int_{\sigma=0}^{\tau}\beta(\sigma)d\sigma\right)[\alpha g_0'(\tau)+\beta(\tau)g_0(\tau)]d\tau$$

$$= \exp\left(\omega(\lambda)\int_{\tau=0}^{t}\beta(\tau)d\tau\right)\int_{\tau=-T}^{0}\exp\left(-\omega(\lambda)\int_{\sigma=0}^{\tau}\beta(\sigma)d\sigma\right)[\alpha g_0'(\tau)+\beta(\tau)g_0(\tau)]d\tau$$

$$= e^{\omega(\lambda)B(t)}\int_{\tau=-T}^{0}e^{-\omega(\lambda)B(\tau)}[\alpha g_0'(\tau)+\beta(\tau)g_0(\tau)]d\tau .$$

Therefore,

$$I_3(x,t+T)-I_3(x,t)=\int_{\mathcal{C}}\left\{\frac{i\lambda}{1+\alpha\lambda^2}e^{i\lambda x+B(t)\omega(\lambda)}\int_{\tau=-T}^{0}e^{-\omega(\lambda)B(\tau)}[\alpha g_0'(\tau)+\beta(\tau)g_0(\tau)]d\tau\right\}d\lambda$$

$$=\int_{\mathcal{C}}e^{\omega(\lambda)B(t)}h(\lambda)d\lambda=\int_{-1/\sqrt{\alpha}}^{1/\sqrt{\alpha}}e^{\omega(\lambda)B(t)}h(\lambda)d\lambda+\int_{\mathcal{K}}e^{\omega(\lambda)B(t)}h(\lambda)d\lambda ,$$

where

$$h(\lambda):=\frac{i\lambda}{1+\alpha\lambda^2}e^{i\lambda x}\int_{\tau=-T}^{0}e^{-\omega(\lambda)B(\tau)}[\alpha g_0'(\tau)+\beta(\tau)g_0(\tau)]d\tau .$$

Computing

$h(0)=0$,

$$h'(\lambda)=\frac{i(1-\alpha\lambda^2)}{(1+\alpha\lambda^2)^2}e^{i\lambda x}\int_{\tau=-T}^{0}e^{-\omega(\lambda)B(\tau)}[\alpha g_0'(\tau)+\beta(\tau)g_0(\tau)]d\tau$$

$$+\frac{i\lambda}{1+\alpha\lambda^2}ixe^{i\lambda x}\int_{\tau=-T}^{0}e^{-\omega(\lambda)B(\tau)}[\alpha g_0'(\tau)+\beta(\tau)g_0(\tau)]d\tau$$

$$+\frac{2i\lambda^2B(\tau)}{(1+\alpha\lambda^2)^3}e^{i\lambda x}\int_{\tau=-T}^{0}e^{-\omega(\lambda)B(\tau)}[\alpha g_0'(\tau)+\beta(\tau)g_0(\tau)]d\tau$$

and

$$h''(0)=-2x\int_{-T}^{0}\beta(\tau)g_0(\tau)d\tau ,$$

we find

$$-2[I_3(x,t+T)-I_3(x,t)]=\left(\sqrt{\pi}x\int_{-T}^{0}\beta(\tau)g_0(\tau)d\tau\right)\frac{1}{[B(t)]^{3/2}}+\mathrm{O}(1/[B(t)]^{5/2}). \tag{5.15}$$

*Step 3* We have

$$I_4(x,t)-I_5(x,t)$$

$$=\int_{-\infty}^{\infty}\left\{e^{i\lambda x+B(t)\omega(\lambda)}\int_{\tau=0}^{t}\frac{e^{-B(\tau)\omega(\lambda)}}{1+\alpha\lambda^2}\hat{f}(\lambda,\tau)\,d\tau\right\}d\lambda-\int_{\mathcal{C}}\left\{e^{i\lambda x+B(t)\omega(\lambda)}\int_{\tau=0}^{t}\frac{e^{-B(\tau)\omega(\lambda)}}{1+\alpha\lambda^2}\hat{f}(-\lambda,\tau)\,d\tau\right\}d\lambda$$

and, by the periodicity assumptions,

$$[I_4(x,t+T)-I_5(x,t+T)]-[I_4(x,t)-I_5(x,t)]$$

$$=\int_{-\infty}^{\infty}e^{B(t)\omega(\lambda)}\left\{\frac{e^{i\lambda x}}{1+\alpha\lambda^2}\int_{\tau=-T}^{0}e^{-B(\tau)\omega(\lambda)}\hat{f}(\lambda,\tau)\,d\tau\right\}d\lambda-\int_{\mathcal{C}}e^{B(t)\omega(\lambda)}\left\{\frac{e^{i\lambda x}}{1+\alpha\lambda^2}\int_{\tau=-T}^{0}e^{-B(\tau)\omega(\lambda)}\hat{f}(-\lambda,\tau)\,d\tau\right\}d\lambda .$$

In view of the above formula, working as in steps 1 and 2, we show that

$$[I_4(x,t+T)-I_5(x,t+T)]-[I_4(x,t)-I_5(x,t)]$$
$$=\int_{-1/\sqrt{\alpha}}^{1/\sqrt{\alpha}} e^{B(t)\omega(\lambda)}\left\{\frac{e^{i\lambda x}}{1+\alpha\lambda^2}\int_{\tau=-T}^{0} e^{-B(\tau)\omega(\lambda)}[\hat{f}(\lambda,\tau)-\hat{f}(-\lambda,\tau)]d\tau\right\}d\lambda+\mathrm{O}(1/[B(t)]^{\ell})$$
$$=\int_{-1/\sqrt{\alpha}}^{1/\sqrt{\alpha}} e^{B(t)\omega(\lambda)}h(\lambda)d\lambda+\mathrm{O}(1/[B(t)]^{\ell}),\ \text{for every } \ell\in\mathbb{N},$$

where we have set – this time –

$$h(\lambda):=\frac{e^{i\lambda x}}{1+\alpha\lambda^2}\int_{\tau=-T}^{0} e^{-B(\tau)\omega(\lambda)}[\hat{f}(\lambda,\tau)-\hat{f}(-\lambda,\tau)]d\tau$$
$$=\frac{e^{i\lambda x}}{1+\alpha\lambda^2}\int_{\tau=-T}^{0} e^{-B(\tau)\omega(\lambda)}\left[\int_0^{\infty} e^{-i\lambda y}f(y,\tau)dy-\int_0^{\infty} e^{i\lambda y}f(y,\tau)dy\right]d\tau.$$

Now we compute:

$$h(0)=0,$$

$$h'(\lambda)=ix\frac{e^{i\lambda x}}{1+\alpha\lambda^2}\int_{\tau=-T}^{0} e^{-B(\tau)\omega(\lambda)}[\hat{f}(\lambda,\tau)-\hat{f}(-\lambda,\tau)]d\tau$$
$$-2\alpha\lambda\frac{e^{i\lambda x}}{(1+\alpha\lambda^2)^2}\int_{\tau=-T}^{0} e^{-B(\tau)\omega(\lambda)}[\hat{f}(\lambda,\tau)-\hat{f}(-\lambda,\tau)]d\tau$$
$$+\frac{2\lambda e^{i\lambda x}}{(1+\alpha\lambda^2)^3}\int_{\tau=-T}^{0} B(\tau)e^{-B(\tau)\omega(\lambda)}[\hat{f}(\lambda,\tau)-\hat{f}(-\lambda,\tau)]d\tau$$
$$+\frac{e^{i\lambda x}}{1+\alpha\lambda^2}\int_{\tau=-T}^{0} e^{-B(\tau)\omega(\lambda)}\left[-i\int_0^{\infty} ye^{-i\lambda y}f(y,\tau)dy-i\int_0^{\infty} ye^{i\lambda y}f(y,\tau)dy\right]d\tau$$

and

$$h''(0)=4x\int_{-T}^{0}\int_0^{\infty} yf(y,\tau)dyd\tau.$$

Thus, again from Laplace formula, we obtain

$$(5.16)\quad [I_4(x,t+T)-I_5(x,t+T)]-[I_4(x,t)-I_5(x,t)]$$
$$=\left(\sqrt{\pi}x\int_{-T}^{0}\int_0^{\infty} yf(y,\tau)dyd\tau\right)\frac{1}{[B(t)]^{3/2}}+\mathrm{O}(1/[B(t)]^{5/2}).$$

*Step 4* The conclusion of the theorem follows from formulas (5.14), (5.15) and (5.16).

*Remark* In view of the remark following the statement of the theorem, formula (5.10) can be written also as follows:

$$u(x,t+T)-u(x,t)=\frac{x}{2\mu^{3/2}\sqrt{\pi}}\left(\int_{-T}^{0}\beta(\tau)g_0(\tau)d\tau+\int_{-T}^{0}\int_0^{\infty} yf(y,\tau)dyd\tau\right)\frac{1}{t^{3/2}}+\mathrm{O}(1/t^{5/2}),$$

where

$$\mu=\frac{1}{T}\int_0^{T}\beta(\tau)d\tau.$$

*Remark* For rigorous UTM-based eventual-periodicity results in autonomous cases see [30,28,18,24].

*Acknowledgment* AC and TO are supported by JSPS KAKENHI Grant Number 24H00024 and JST [Moonshot R&D][Grant Number JPMJMS25A5] .

## References

[1] A. S. Fokas, A unified transform method for solving linear and certain nonlinear PDEs, Proc. Roy. Soc. London A, 453 (1997) 1411–1443.
[2] A.S. Fokas, "On the integrability of certain linear and nonlinear partial differential equations", *J. Math. Phys.* 41, 4188–4237 (2000);
[3] A.S. Fokas, "A new transform method for evolution partial differential equations", *IMA J. Appl. Math.* 67, 559-590 (2002).
[4] Fokas, A. S. (2004). Boundary-value problems for linear PDEs with variable coefficients. *Proceedings of the Royal Society of London. Series A: Mathematical, Physical and Engineering Sciences*, *460*(2044), 1131-1151.
[5] A. S. Fokas, "From Green to Lax via Fourier", In: *Recent Advances in Nonlinear Partial Differential Equations and Applications*, *Proc. Sympos. Appl. Math.* 65 (Dedicated to P.D. Lax and L. Nirenberg), Amer. Math. Soc. (2007).
[6] A. S. Fokas, A Unified Approach to Boundary Value Problems, (CBMS-NSF Regional Conference Series in Applied Mathematics) Philadelphia, SIAM; (2008).
[7] A.S. Fokas, "Lax Pairs: A novel type of separability" (invited paper), *Inverse Problems* 25 (2009).
[8] A. S. Fokas, E. A. Spence, Synthesis, as opposed to separation, of variables, *SIAM Review* 54 (2012).
[9] A Chatziafratis, AS Fokas, K Kalimeris, The Fokas method for evolution partial differential equations, Partial Differ. Equ. Appl. Math. (2025).
[10] Chatziafratis, A., Fokas, A. S., & Aifantis, E. C. (2024). Variations of heat equation on the half-line via the Fokas method (invited review article). *Mathematical Methods in the Applied Sciences*, *47*(16), 12541-12566.
[11] J. L. Bona, A. Chatziafratis, H. Chen, S. Kamvissis, The linearised BBM equation on the half-line revisited, *Lett. Math. Phys.* 114 (2024).
[12] M. M. Cavalcanti, A. Chatziafratis, and L. Rosier. A unified approach to lack of controllability for linear evolution equations on quarter-planes. *preprint* (2025), to appear.
[13] Chatziafratis, A. Fokas type closed-form solution formulae for Sobolev-type equations with time-dependent coefficients. *arXiv preprint arXiv:2512.15937* (2025).
[14] A. Chatziafratis et al., Rigorous analysis for the Dirac system on the quarter-plane, preprint, 2026.
[15] Chatziafratis, A., & Aifantis, E. C. (2025). Explicit Ehrenpreis–Palamodov–Fokas representations for the Sobolev–Barenblatt model and the Rubinshtein–Aifantis system on semi-strips and quarter-planes. *Lobachevskii Journal of Mathematics*, *46*(2), 594-612.
[16] A. Chatziafratis, E.C. Aifantis, A. Carbery, A.S. Fokas, Integral representations for the double-diffusion system on the half-line, *Z. Angew. Math. Phys.* 75 (2024).
[17] Chatziafratis, A., Fokas, A.S., Aifantis, E.C. (2024). On Barenblatt's pseudoparabolic equation with forcing on the half-line via the Fokas method. *Z. Angew. Math. Mech. 104*(3), e202300614.
[18] Chatziafratis, A., Giorgi, C., Miranville, A., & Zullo, F. (2026). On generalised d'Alembert-type integral representations for damped wave equations on the quarter-plane, *arXiv:2608.06355*
[19] A. Chatziafratis, L. Grafakos, S. Kamvissis, Long-range instabilities for linear evolution PDE on semi-bounded domains, *Dyn. PDE* 21 (2024).
[20] Chatziafratis, A., & Hwang, G., On a one-phase Stefan-type free boundary problem with advection and reaction, *Computational and Applied Mathematics* (2026).
[21] A. Chatziafratis, S. Kamvissis, I. G. Stratis, Boundary behavior of the solution to the linear KdV equation on the half-line, Stud. Appl. Math. (2023).
[22] Chatziafratis, A., & Kamvissis, S. (2024). A note on uniqueness for linear evolution PDEs posed on the quarter-plane. *arXiv preprint arXiv:2401.08531*.
[23] A. Chatziafratis, S. Kamvissis, Infinity of solutions to initial-boundary value problems for linear constant-coefficient evolution PDEs on semi-infinite intervals, *Bull. London Math. Soc.* (2025).

[24] Chatziafratis, A., Kamvissis, S., & Slavík, A. (2026). Initial-boundary value problems for linear lattice differential equations. *arXiv preprint arXiv:2610.31057*.
[25] A. Chatziafratis, G. Karali, C. E. Synolakis (2026). Rigorous Analytical Solution for the Linearized Whitham–Broer–Kaup System on the Half-Line. *Studies in Applied Mathematics*, *156*(6).
[26] A. Chatziafratis, E. Luca and S. L. Waters, The Milne–Barenblatt pseudo-parabolic equation with time-dependent coefficients in biological fluid mechanics and regenerative medicine, *Discrete Contin Dyn Syst-S* (2026).
[27] Chatziafratis, A., & Mantzavinos, D. (2022). Boundary behavior for the heat equation on the half-line. *Mathematical Methods in the Applied Sciences*, *45*(12), 7364-7393.
[28] A. Chatziafratis, A. Miranville, G. Karali, A.S. Fokas, E.C. Aifantis, Higher-order diffusion and Cahn–Hilliard-type models revisited on the half-line, *Math. Models Methods Appl. Sci.* (2025).
[29] Chatziafratis, A., Miranville, A., & Ozawa, T. (2026). The linear Cahn-Hilliard equation with an interface. *arXiv preprint arXiv:2605.20130*.
[30] A. Chatziafratis, T. Ozawa, New instability, blow-up and break-down effects for Sobolev-type evolution PDE: asymptotic analysis for a celebrated pseudo-parabolic model on the quarter-plane. *Partial Differential Equations and Applications*, *5*(5), 30 (2024).
[31] A. Chatziafratis, T. Ozawa, Non-standard initial-boundary-value and interface problems for the linear BBM equation, *Discrete Contin Dyn Syst-S* (2026).
[32] A. Chatziafratis, T. Ozawa, S.-F. Tian, Rigorous analysis of the unified transform method and long-range instabilities for the inhomogeneous time-dependent Schrödinger equation on the quarter-plane, *Math. Annalen* 389 (2023).
[33] Chatziafratis, A., Rukolaine, S. A., & Aifantis, E. C. (2026). Integral representations for the Barenblatt–Sobolev–Galpern pseudo-parabolic equation: boundary-value and interface problems on unbounded and finite intervals. *Russ J Math Phys*, *33*(3), 465-491.
[34] Milne, E. A., The diffusion of imprisoned radiation through a gas. *J. London Math. Soc*. (1926)
[35] Barenblatt, G. I., Zheltov, Iu. P., Kochina, I. N. (1960) Basic concepts in the theory of seepage of homogeneous liquids in fissured rocks (strata). *J. Appl. Math. Mech.*, 24: 1286–1303.
[36] Ting, T. W., Certain non-steady flows of second-order fluids. Arch. Rational Mech. Anal. (1963)
[37] Showalter, R. E., Ting, T. W. (1970) Pseudoparabolic partial differential equations. *SIAM J. Math. Anal.*, 1: 1–26.
[38] Puri, P., Impulsive motion of a flat plate in a Rivlin-Ericksen fluid. *Rheol. Acta*, 23 (1984).
[39] Chen, P. J., Gurtin, M. E. (1968) On a theory of heat conduction involving two temperatures. *Z. Angew. Math. Phys.*, 19: 614–627.
[40] D.E. Amos, On half-space solutions of a modified heat equation, *Quart Appl Math* 28 (1969).
[41] R. E. Showalter, T.W. Ting, Asymptotic behaviour of solutions of pseudo-parabolic partial differential equations, *Ann Mat Pura Appl* 90, 241-258 (1971).
[42] D. Colton, Pseudoparabolic equations in one space variable, J Differ Equ 12(3), 559-565 (1972).
[43] V. R. Gopala Rao, T.W. Ting, Solutions of pseudo-heat equations in the whole space, *Arch Rat Mech Anal* 49 (1972).
[44] T. W. Ting, A cooling process according to two-temperature theory of heat conduction, *J Math Anal Appl* 45, 23-31 (1974).
[45] W. Rundell, The uniqueness class for the Cauchy problem for pseudoparabolic equations, *Proc Am Math Soc* 76 (1979).
[46] E. di Benedetto, M. Pierre, On the maximum principle for pseudoparabolic equations, *Indiana Univ Math J* 30, 821 (1980).
[47] L.W. White, Controllability properties of pseudo-parabolic boundary control problems, *SIAM J Control Optimiz* 18 (1980).
[48] J. R. Cannon, Y. Lin, Classical and weak solutions for one-dimensional pseudo-parabolic equations with typical boundary data, *Ann Mat Pura Appl* 152, 375-385 (1987).
[49] Chatziafratis, A., & Hatziafratis, T., The natural boundary of Fourier-Laplace transforms of functions in the half-line Schwartz space. *Integral Transforms and Special Functions* (2025).